\documentclass[a4paper, 12pt]{article}%
\usepackage{amsmath}
\usepackage{graphicx}%
\usepackage{amsfonts}%
\usepackage{amssymb}
\begin{document}

\begin{center}
{\Huge Bergelson's Theorem for weakly mixing C*-dynamical systems}

\bigskip

Rocco Duvenhage

\bigskip

\textit{Department of Mathematics and Applied Mathematics}

\textit{University of Pretoria, 0002 Pretoria, South Africa}

\bigskip

2006-7-13
\end{center}

\bigskip

\noindent\textbf{Abstract}

We study a nonconventional ergodic average for asymptotically abelian weakly
mixing C*-dynamical systems, related to a second iteration of Khintchine's
recurrence theorem obtained by Bergelson in the measure theoretic case. A
noncommutative recurrence theorem for such systems is obtained as a corollary.

\bigskip


\section{Introduction}

In 1977 Furstenberg \cite{F77} published a very influential paper where he
proved a recurrence theorem for measure preserving dynamical systems $\left(
X,\Sigma,\nu,T\right)  $, which follows from%
\begin{equation}
\liminf_{N\rightarrow\infty}\frac{1}{N}\sum_{n=1}^{N}\nu\left(  V\cap
T^{-n}V\cap T^{-2n}V\cap...\cap T^{-kn}V\right)  >0 \tag{1.1}%
\end{equation}
for any $V\in\Sigma$ with $\nu(V)>0$, and lead to an alternative proof of
Szemer\'{e}di's Theorem in combinatorial number theory. This approach to
Szemer\'{e}di's Theorem lead to various generalizations of the latter and a
field of research now often called Ergodic Ramsey Theory. Recently Niculescu,
Str\"{o}h and Zsid\'{o} \cite{NSZ} initiated\ a programme to extend
Furstenberg's result to C*-dynamical systems, and more generally to study
``noncommutative recurrence''; also see \cite{NS} and \cite{deBDS}.

Meanwhile much research, e.g. \cite{F90, BHK, HK05}, has been done for measure
theoretic dynamical systems to determine when the $\lim\inf$ in (1.1), and
generalizations thereof, is in fact a limit (the study of ``nonconventional
ergodic averages''), and to find lower bounds for these limits similar to the
lower bound appearing in
\begin{equation}
\lim_{N\rightarrow\infty}\frac{1}{N}\sum_{n=1}^{N}\nu\left(  V\cap
T^{-n}V\right)  \geq\nu(V)^{2} \tag{1.2}%
\end{equation}
which follows from the mean ergodic theorem, and from which in turn
Khintchine's recurrence theorem follows. In particular in \cite{HK05} it was
shown that (1.1) is indeed a limit, but certain negative results regarding
lower bounds were found in \cite{BHK}.

However a very interesting theorem was proven by Bergelson \cite{B}, which was
later significantly generalized by Host and Kra \cite{HK04, HK05}. Simply put
they consider averages along cubes in $\mathbb{Z}^{q}$, rather than along
arithmetic progressions as in (1.1). In particular Bergelson's Theorem covers
the two dimensional case, i.e. a square in $\mathbb{Z}^{2}$, which can also be
viewed as a second iteration of (1.2), with the average being of the form%
\begin{align}
&  \lim_{N\rightarrow\infty}\frac{1}{N^{2}}\sum_{m=1}^{N}\sum_{n=1}^{N}%
\nu\left(  V\cap T^{-n}V\cap T^{-m}\left(  V\cap T^{-n}V\right)  \right)
\nonumber\\
&  =\lim_{N\rightarrow\infty}\frac{1}{N^{2}}\sum_{m=1}^{N}\sum_{n=1}^{N}%
\nu\left(  V\cap T^{-n}V\cap T^{-m}V\cap T^{-(m+n)}V\right) \tag{1.3}\\
&  \geq\nu(V)^{4}\text{.}\nonumber
\end{align}

For simplicity the averages in (1.2) and (1.3) were taken over $\left[
1,N\right]  $ and $\left[  1,N\right]  \times\left[  1,N\right]  $
respectively. But in fact, the average in (1.2) can be taken over $\left[
M,N\right]  $, and as \cite{B} shows, the average in (1.3) can be taken over
$\left[  M,N\right]  \times\left[  M,N\right]  $, with the limit
$N-M\rightarrow\infty$ being taken, and the results then still hold. This
provides a uniformness which leads to the relative denseness in the resulting
recurrence theorems, for example in Khintchine's case for any $\varepsilon>0$
the set
\[
\left\{  n:\nu\left(  V\cap T^{-n}V\right)  >\nu(V)^{2}-\varepsilon\right\}
\]
is relatively dense (also said to be syndetic) in $\mathbb{N}=\{1,2,3,...\}$,
i.e. the set has bounded gaps.

In this paper we study an extension of Bergelson's Theorem to C*-dynamical
systems. The main difference of course is that the probability space $\left(
X,\Sigma,\nu\right)  $ and some abelian algebra of functions on it, like
$L^{\infty}(\nu)$, are replaced by a unital C*-algebra $A$ which need not be
abelian, and a state $\omega$ on $A$. We also work over more general abelian
groups than $\mathbb{Z}$. We follow the basic structure of Bergelson's proof
\cite[Section 5]{B}. However we will only prove Bergelson's Theorem for
asymptotically abelian (in the sense of Definition 2.5) weakly mixing
C*-dynamical systems with $\omega$ a trace, so a degree of abelianness is
still present in the dynamical system. In order to get the uniformness
mentioned above, we need to restrict further to countable groups, and use a
stronger form of asymptotic abelianness which we call ``uniform asymptotic
abelianness''; see Definition 2.5.

The main results are Theorem 5.6 and its corollary at the end of the paper.
The rest of the paper systematically builds up the required tools, including
an appropriate ``van der Corput lemma'' in Section 5, to prove this theorem.

\section{Definitions and notations}

In this section we collect some of the definitions and notations that we will
use in the rest of the paper.

To simplify statements of definitions and results, we introduce the following
terminology and notation:

The space of all linear mappings $X\rightarrow X$ on a vector space $X$ will
be denoted by $L(X)$, while the space of all bounded linear operators
$X\rightarrow X$ on a normed space $X$ will be denoted by $B(X)$. For a $\ast
$-algebra $A$ we will denote the group of all $\ast$-isomorphisms
$A\rightarrow A$ by Aut$(A)$.

For a state $\omega$ on a unital $\ast$-algebra $A$, i.e. a linear functional
on $A$ such that $\omega(a^{\ast}a)\geq0$ and $\omega(1)=1$, we will denote
any GNS representation of $(A,\omega)$ by $(H,\iota)$, which is a Hilbert
space $H$ and a linear mapping $\iota:A\rightarrow H$ such that $\left\langle
\iota(a),\iota(b)\right\rangle =\omega(a^{\ast}b)$ for all $a,b\in A$ and with
$\iota(A)$ dense in $H$. Note that we use the convention where inner products
are conjugate linear in the first slot. This mapping $\iota$ can be expressed
in terms of a homomorphism $\pi:A\rightarrow L(\iota(A))$, or by a $\ast
$-homomorphism $\pi:A\rightarrow B(H)$ when $A$ is a C*-algebra, and the
formula $\iota(a)=\pi(a)\Omega$ where $\Omega:=\iota(1)\in H$. More generally
$\pi$ is given by $\pi(a)\iota(b)=\iota(ab)$. Related to this, we note that
$\left|  \left|  a\right|  \right|  _{\omega}:=\sqrt{\omega\left(  a^{\ast
}a\right)  }=\left|  \left|  \iota(a)\right|  \right|  $ defines a seminorm on
$A$, and when $A$ is a C*-algebra we have $\left|  \left|  a\right|  \right|
_{\omega}\leq\left|  \left|  a\right|  \right|  $. A state $\omega$ will be
called \textit{tracial} or a \textit{trace} when $\omega(ab)=\omega(ba)$ for
all $a,b\in A$.

In a group or semigroup $G$ we will use the notations $Vg:=\left\{  vg:v\in
V\right\}  $, $VW:=\left\{  vw:v\in V,w\in W\right\}  $, $V^{-1}:=\left\{
v^{-1}:v\in V\right\}  $ (in groups), etc. for any $V,W\subset G$ and $g\in
G$, and we will use multiplicative notation even when the group or semigroup
is abelian. If $\Sigma$ is a $\sigma$-algebra in some set, then for any
$V\in\Sigma$ we have a $\sigma$-algebra $\Sigma|_{V}:=\left\{  W\cap
V:W\in\Sigma\right\}  $ in $V$, and we let $\Sigma\times\Sigma$ or $\Sigma
^{2}$ denote the product $\sigma$-algebra of $\Sigma$ with itself. A triple
$\left(  K,\Sigma,\mu\right)  $ with $K$ a semigroup, $\Sigma$ a $\sigma
$-algebra in $K$ with $\Sigma g\subset\Sigma$ (equality being allowed by this
notation) for all $g\in K$, and $\mu$ a positive measure on $\Sigma$, will be
called a \textit{measure semigroup}. When $\mu$ is right invariant, i.e.
$\mu(Vg)=\mu(V)$ for $V\in\Sigma$ and $g\in K$, we say that $\left(
K,\Sigma,\mu\right)  $ is \textit{right invariant}, and similarly for
\textit{left invariance}. A measure semigroup with the semigroup being a
group, will simply be called a \textit{measure group}. In integrals with
respect to the given measure $\mu$ we will often write $dg$ with $g$ being the
variable involved, instead of $d\mu$ or $d\mu(g)$.

When we say that a net $\left(  \Lambda_{\alpha}\right)  $ has some property
for $\alpha$ ``large enough'', then we mean that there is a $\beta$ in the
directed set such that the property holds for all $\alpha>\beta$.

\bigskip

\noindent\textbf{Definition 2.1.} Let $\left(  K,\Sigma,\mu\right)  $ be a
measure semigroup. A net $\left(  \Lambda_{\alpha}\right)  $ in $\Sigma$ is
called a \textit{F\o lner net} (or a \textit{F\o lner sequence}, in case the
directed set is $\mathbb{N}$) in $\left(  K,\Sigma,\mu\right)  $ if
$0<\mu(\Lambda_{\alpha})<\infty$ for $\alpha$ large enough and
\[
\lim_{\alpha}\frac{\mu\left(  \Lambda_{\alpha}\Delta(\Lambda_{\alpha
}g)\right)  }{\mu(\Lambda_{\alpha})}=0
\]
for all $g\in K$. (In some cases $\Lambda_{\alpha}$ will be constructed from
another F\o lner net, and in such cases $\Lambda_{\alpha}$ will only be
required to be in $\Sigma$ for $\alpha$ large enough.) A F\o lner net $\left(
\Lambda_{\alpha}\right)  $ is called \textit{uniform} if%
\[
\lim_{\alpha}\frac{1}{\mu(\Lambda_{\alpha})}\sup_{g\in\Lambda_{\beta}}%
\mu\left(  \Lambda_{\alpha}\Delta(\Lambda_{\alpha}g)\right)  =0
\]
for $\beta$ large enough.

\bigskip

\noindent\textbf{Definition 2.2.} Let $\left(  G,\Sigma,\mu\right)  $ be a
right invariant measure group. Let $K\in\Sigma$ be a semigroup, and assume
that $\left(  K,\Sigma|_{K},\mu|_{\Sigma|_{K}}\right)  $ has a F\o lner net
$\left(  \Lambda_{\alpha}\right)  $, then we will call $K$ a \textit{F\o lner
semigroup}. In particular, when $K=G$, we call $G$ a \textit{F\o lner group}.
(In the sequel, whenever we refer to a F\o lner semigroup $K$, such an
``ambient'' group $\left(  G,\Sigma,\mu\right)  $ as well as $K$ 's measure
space structure $\left(  K,\Sigma|_{K},\mu|_{\Sigma|_{K}}\right)  $ will be
implied and in such cases we will use this notation consistently. When we
refer to a F\o lner group $G$, then it is its own ambient group.)

\bigskip

\noindent\textbf{Definition 2.3.} A F\o lner net $\left(  \Lambda_{\alpha
}\right)  $ in a F\o lner semigroup $K$ is said to satisfy the
\textit{Tempel'man condition} if there is a real number $c>0$ such that
\[
\mu\left(  \Lambda_{\alpha}^{-1}\Lambda_{\alpha}\right)  \leq c\mu\left(
\Lambda_{\alpha}\right)
\]
for $\alpha$ large enough (in particular $\Lambda_{\alpha}^{-1}\Lambda
_{\alpha}\in\Sigma$ is required for $\alpha$ large enough).

\bigskip

Now we define the types of dynamical systems with which we will be working:

\bigskip

\noindent\textbf{Definition 2.4. }Let $\omega$ be a state on a unital $\ast
$-algebra $A$. Let $\left(  K,\Sigma,\mu\right)  $ be a measure semigroup and
consider a function $\tau:K\rightarrow L(A):g\mapsto\tau_{g}$ with
\begin{align*}
\tau_{g}\circ\tau_{h}  &  =\tau_{gh}\\
\tau_{g}(1)  &  =1\\
\left\|  \tau_{g}(a)\right\|  _{\omega}  &  \leq\left\|  a\right\|  _{\omega}%
\end{align*}
for all $g,h\in K$ and $a\in A$, and such that $K\rightarrow\mathbb{C}%
:g\mapsto\omega\left(  a\tau_{g}(b)\right)  $ is $\Sigma$-measurable for all
$a,b\in A$. Then we'll call $\left(  A,\omega,\tau,K\right)  $ a $\ast
$\textit{-dynamical system}, the structure $\left(  \Sigma,\mu\right)  $ on
$K$ being implied by this notation. The following special cases will be
important (note that here we use arbitrary nets in $\Sigma$ rather than F\o
lner nets, and also note that (2) implies (1)):

(1) If
\[
\lim_{\alpha}\frac{1}{\mu\left(  \Lambda_{\alpha}\right)  }\int_{\Lambda
_{\alpha}}\omega\left(  a\tau_{g}(b)\right)  dg=\omega(a)\omega(b)
\]
for all $a,b\in A$ for some net $\left(  \Lambda_{\alpha}\right)  $ in
$\Sigma$ with $0<\mu(\Lambda_{\alpha})<\infty$ for $\alpha$ large enough, then
we will call $\left(  A,\omega,\tau,K\right)  $ \textit{ergodic with respect
to }$\left(  \Lambda_{\alpha}\right)  $.

(2) If
\[
\lim_{\alpha}\frac{1}{\mu\left(  \Lambda_{\alpha}\right)  }\int_{\Lambda
_{\alpha}}\left|  \omega\left(  a\tau_{g}(b)\right)  -\omega(a)\omega
(b)\right|  dg=0
\]
for all $a,b\in A$ for some net $\left(  \Lambda_{\alpha}\right)  $ in
$\Sigma$ with $0<\mu(\Lambda_{\alpha})<\infty$ for $\alpha$ large enough, then
we will call $\left(  A,\omega,\tau,K\right)  $ \textit{weakly mixing with
respect to }$\left(  \Lambda_{\alpha}\right)  $.

\bigskip

When working with such systems, we can use any GNS representation to represent
$\tau$ on $H$ by the formula
\[
U_{g}\iota(a):=\iota\left(  \tau_{g}(a)\right)
\]
and then uniquely extending $U_{g}:\iota(A)\rightarrow\iota(A)$ to $H$ to give
us a \textit{representation} $U:K\rightarrow B(H):g\mapsto U_{g}$ of the
semigroup $K$ as contractions, i.e. $U_{g}U_{h}=U_{gh}$ and $\left|  \left|
U_{g}\right|  \right|  \leq1$ for all $g,h\in K$. We'll consider this $U$ to
be the corresponding \textit{GNS representation} of $\tau$, and in the
presence of a GNS representation we will use this notation in the rest of the
paper. Note that $U_{g}\Omega=\iota\left(  \tau_{g}(1)\right)  =\iota
(1)=\Omega$.

\bigskip

\noindent\textbf{Definition 2.5.} Let $\omega$ be a state on a unital
C*-algebra $A$ and let $K$ be a F\o lner semigroup with $G$ abelian. Assume
that $\mathcal{I}:G\rightarrow G:g\mapsto g^{-1}$ is measurable, i.e.
$\mathcal{I}^{-1}(\Sigma)\subset\Sigma$, and $\mu$-invariant, i.e. $\mu
\circ\mathcal{I}=\mu$. Denote the identity of $G$ by $e$. Consider a function
$\tau:G\rightarrow$Aut$(A):g\mapsto\tau_{g}$ such that $\tau_{e}$ is the
identity on $A$, $\tau_{g}\circ\tau_{h}=\tau_{gh}$ and $\omega\circ\tau
_{g}=\omega$ for all $g,h\in G$, and such that $K\rightarrow\mathbb{C}%
:g\mapsto\omega\left(  a\tau_{g}(b)\right)  $ is $\left(  \Sigma|_{K}\right)
$-measurable for all $a,b\in A$. Then we'll call $\left(  A,\omega
,\tau,K\right)  $ a \textit{C*-dynamical system}. (It is easily seen that this
is a special case of a $\ast$-dynamical system.) If furthermore $K\rightarrow
\mathbb{R}:g\mapsto\left\|  \left[  a,\tau_{g}(b)\right]  \right\|  $ is
$\left(  \Sigma|_{K}\right)  $-measurable for all $a,b\in A$, where $\left[
a,b\right]  :=ab-ba$, and $\left(  \Lambda_{\alpha}\right)  $ is a F\o lner
net in $K$, then the following special cases will be considered:

(1) We call $\left(  A,\omega,\tau,K\right)  $ \textit{asymptotically abelian
with respect to }$\left(  \Lambda_{\alpha}\right)  $ if
\[
\lim_{\alpha}\frac{1}{\mu\left(  \Lambda_{\alpha}\right)  }\int_{\Lambda
_{\alpha}}\left\|  \left[  a,\tau_{g}(b)\right]  \right\|  dg=0
\]
for all $a,b\in A$.

(2) We call $\left(  A,\omega,\tau,K\right)  $ \textit{uniformly
asymptotically abelian with respect to }$\left(  \Lambda_{\alpha}\right)  $
if
\[
\lim_{\alpha}\sup_{h\in K}\frac{1}{\mu\left(  \Lambda_{\alpha}\right)  }%
\int_{\Lambda_{\alpha}h}\left\|  \left[  a,\tau_{g}(b)\right]  \right\|  dg=0
\]
for all $a,b\in A$.

\bigskip

\noindent\textbf{Remarks.} The terminology in Definition 2.5(1) and (2) isn't
quite standard, but we'll use it consistently in this paper. For simplicity,
in Definition 2.5 consider the case $K=G=\mathbb{Z}$ and the F\o lner sequence
in $\mathbb{Z}$ given by $\Lambda_{N}=\left\{  1,...,N\right\}  $. The term
``asymptotic abelian'' (see for example \cite{DKS}) is often used to describe
the condition
\[
\lim_{\left|  n\right|  \rightarrow\infty}\left\|  \left[  a,\tau
_{n}(b)\right]  \right\|  =0
\]
which for the purposes of these remarks we'll refer to as ``strong asymptotic
abelianness''. Note that this condition implies asymptotic abelianness and
uniform asymptotic abelianness in our sense above with respect to $\left(
\Lambda_{N}\right)  $. In fact, if $\omega$ is a so-called factor state, then
strong asymptotic abelianness also implies what is known as ``strong mixing''
namely
\[
\lim_{\left|  n\right|  \rightarrow\infty}\left|  \omega\left(  a\tau
_{n}(b)\right)  -\omega(a)\omega(b)\right|  =0
\]
(see \cite[Example 4.3.24]{BR} for details) and hence weak mixing with respect
to $\left(  \Lambda_{N}\right)  $. We can also mention that from the results
and discussions in \cite[Section 5.3.2]{BR2} regarding infinite temperature
KMS states (i.e. at inverse temperature $\beta=0$), it follows that if a
unital C*-algebra $A$ has at least one trace, then it also has a ``factor
trace'' (i.e. a factor state which is tracial).

\section{Background}

We now discuss a number of results which for the most part are known, but we
formulate and adapt them in a way that will suit our needs in Section 5. Some
more notation is also introduced.

\bigskip

\noindent\textbf{Definition 3.1.} Let $S$ be an abelian semigroup of linear
contractions on a Hilbert space $H$. A vector $x\in H\backslash\left\{
0\right\}  $ is called an \textit{eigenvector with unimodular eigenvalues} (or
``unimodular eigenvector'' for short) of $S$ if there exists a function
$\lambda:S\rightarrow\mathbb{C}$ such that $\left|  \lambda(U)\right|  =1$ and
$Ux=\lambda(U)x$ for every $U\in S$. The closure of the span (finite linear
combinations) of the unimodular eigenvectors of $S$ will be denoted by $H_{0}%
$, and its orthogonal complement in $H$ by $H_{v}$, hence $H=H_{0}\oplus
H_{v}$. The elements of $H_{0}$ are called \textit{reversible}, while the
elements of $H_{v}$ are called \textit{flight vectors}.

\bigskip

An important characterization of $H_{v}$ in our setting is the following:

\bigskip

\noindent\textbf{Proposition 3.2.} \textit{Let }$K$\textit{ be an abelian F\o
lner\ semigroup (though }$G$\textit{ need not be abelian) and }$\left(
\Lambda_{\alpha}\right)  $\textit{ any F\o lner net in }$K$\textit{. Consider
a representation }$U:K\rightarrow B(H)$\textit{ of }$K$\textit{ as
contractions on any Hilbert space }$H$\textit{, giving the semigroup
}$S=\left\{  U_{g}:g\in K\right\}  $\textit{ in Definition 3.1, such that
}$g\mapsto\left\langle x,U_{g}y\right\rangle $\textit{ is }$\left(
\Sigma|_{\Lambda_{\alpha}}\right)  $\textit{-measurable (on }$\Lambda_{\alpha
}$\textit{) and }$\left(  \Sigma|_{\Lambda_{\alpha}h}\right)  $%
\textit{-measurable (on }$\Lambda_{\alpha}h$)\textit{ for all }$h\in
K$\textit{ and }$x,y\in H$\textit{, for }$\alpha$\textit{ large enough. \ Then
for }$y\in H$ \textit{we have the following: }$y\in H_{v}$\textit{ if and only
if }%
\[
\lim_{\alpha}\frac{1}{\mu\left(  \Lambda_{\alpha}\right)  }\int_{\Lambda
_{\alpha}}\left|  \left\langle x,U_{g}y\right\rangle \right|  dg=0
\]
\textit{for all }$x\in H$\textit{.}

\bigskip

This result follows as a special case of the results in \cite[Section 2.4]{K},
but with some modifications to the proof of \cite[Theorem 2.4.7]{K} to
compensate for the fact that we are using the form $\lim_{\alpha}\frac{1}%
{\mu\left(  \Lambda_{\alpha}\right)  }\int_{\Lambda_{\alpha}}(\cdot)d\mu$
rather than an abstract invariant mean as \cite{K} does (the fact that
$\left(  \Lambda_{\alpha}\right)  $ is F\o lner and $\mu$ is invariant, plays
an important role here). An important example is where $U$ is a GNS
representation of a $\ast$-dynamical system, and in such cases we'll use the
notation $H_{0}$ and $H_{v}$ without further explanation. As a corollary we
have the following related characterization of $H_{v}$:

\bigskip

\noindent\textbf{Corollary 3.3.} \textit{Let }$K$\textit{ be an abelian F\o
lner semigroup where }$\left(  G,\Sigma,\mu\right)  $\textit{ is }$\sigma
$\textit{-finite and }$KK=K$\textit{. Consider any F\o lner nets }$\left(
\Lambda_{1\alpha}\right)  $\textit{ and }$\left(  \Lambda_{2\alpha}\right)
$\textit{, indexed by the same directed set, in }$K$\textit{, and a
representation }$U:K\rightarrow B(H)$\textit{ of }$K$\textit{ as contractions
on any Hilbert space }$H$ \textit{such that }$(g,h)\mapsto\left\langle
x,U_{gh}y\right\rangle $\textit{ is }$\left(  \Sigma\times\Sigma\right)
|_{\Lambda_{2\alpha}\times\Lambda_{1\alpha}}$\textit{-measurable and }$\left(
\Sigma\times\Sigma\right)  |_{\Lambda_{2\alpha}h_{2}\times\Lambda_{1\alpha
}h_{1}}$\textit{-measurable for all }$h_{1},h_{2}\in K$\textit{ and }$x,y\in
H$\textit{, for }$\alpha$\textit{ large enough. Then for }$y\in H$ \textit{we
have the following: }$y\in H_{v}$\textit{ if and only if }%
\begin{equation}
\lim_{\alpha}\frac{1}{\mu\left(  \Lambda_{1\alpha}\right)  \mu\left(
\Lambda_{2\alpha}\right)  }\int_{\Lambda_{1\alpha}}\int_{\Lambda_{2\alpha}%
}\left|  \left\langle x,U_{gh}y\right\rangle \right|  dgdh=0 \tag{3.1}%
\end{equation}
\textit{for all }$x\in H$\textit{.}

\bigskip

\noindent\textbf{Proof.} We will apply Proposition 3.2 to $K^{\prime}:=K\times
K$ as a F\o lner semigroup in the ambient group $\left(  G^{\prime}%
,\Sigma^{\prime},\mu^{\prime}\right)  :=\left(  G\times G,\Sigma\times
\Sigma,\mu\times\mu\right)  $ (as in Definition 2.2). It is easily established
that $\Lambda_{\alpha}^{\prime}:=\Lambda_{2\alpha}\times\Lambda_{1\alpha}$
gives a F\o lner net in $K^{\prime}$. That $\left(  G^{\prime},\Sigma^{\prime
},\mu^{\prime}\right)  $ is a right invariant measure group can be shown as
follows: First prove that $\left\{  V\in\Sigma^{\prime}:V(g,h),V^{c}%
(g,h)\in\Sigma^{\prime}\right\}  $ is a $\sigma$-algebra and hence equal to
$\Sigma^{\prime}$ for all $g,h\in G$, where $V^{c}:=G^{\prime}\backslash V$.
From this we obtain $\Sigma^{\prime}(g,h)\subset\Sigma^{\prime}$, and then the
right invariance of $\mu^{\prime}$ follows using the definition of a product measure.

Note that since $K$ is abelian, $U^{\prime}:K^{\prime}\rightarrow B(H)$
defined by $U_{(g,h)}^{\prime}:=U_{gh}$ is a representation of $K^{\prime}$ as
contractions. Since $KK=K$, the semigroups $S:=\left\{  U_{g}:g\in K\right\}
$ and $S^{\prime}:=\left\{  U_{(g,h)}^{\prime}:(g,h)\in K^{\prime}\right\}  $
have the same reversible vectors $H_{0}$, and hence the same flight vectors
$H_{v}$. But then by Proposition 3.2 a $y\in H$ is in $H_{v}$ if and only%
\[
\lim_{\alpha}\frac{1}{\mu^{\prime}\left(  \Lambda_{\alpha}^{\prime}\right)
}\int_{\Lambda_{\alpha}^{\prime}}\left|  \left\langle x,U_{(g,h)}^{\prime
}y\right\rangle \right|  d(g,h)=0
\]
for all $x\in H$, which is exactly (3.1) by Fubini's Theorem. $\square$

\bigskip

We will in fact only need Proposition 3.2 and Corollary 3.3 for measurable
$K\rightarrow\mathbb{C}:g\mapsto\left\langle x,U_{g}y\right\rangle $ and
$K\times K\rightarrow\mathbb{C}:(g,h)\mapsto\left\langle x,U_{gh}%
y\right\rangle $ respectively.

For a F\o lner semigroup $K$ and any F\o lner net $\left(  \Lambda_{\alpha
}\right)  $ in $K$ we have the following mean ergodic theorem: Let $H$ be a
Hilbert space and $U:K\rightarrow B(H):g\mapsto U_{g}$ a representation of $K$
as contractions such that $K\ni g\mapsto\left\langle x,U_{g}y\right\rangle $
is $\left(  \Sigma|_{K}\right)  $-measurable for all $x,y\in H$. Take $P$ to
be the projection of $H$ onto $V:=\left\{  x\in H:U_{g}x=x\text{ for all }g\in
K\right\}  $. Then
\[
\lim_{\alpha}\frac{1}{\mu\left(  \Lambda_{\alpha}\right)  }\int_{\Lambda
_{\alpha}}U_{g}xdg=Px
\]
for all $x\in H$, where integrals over sets $\Lambda\in\Sigma|_{K}$ with
$\mu\left(  \Lambda\right)  <\infty$, of bounded Hilbert space valued
functions $f$ with $\left\langle f,x\right\rangle $ measurable for every $x$
in a dense linear subspace of $H$, are defined via the Riesz representation
theorem, i.e. $\left\langle \int_{\Lambda}fd\mu,x\right\rangle :=\int
_{\Lambda}\left\langle f(g),x\right\rangle dg$. This integral has several
simple properties, for example if $F:K\rightarrow\mathbb{R}$ is measurable and
$\left|  \left|  f\right|  \right|  \leq F$, then $\left|  \left\langle
\int_{\Lambda}fd\mu,x\right\rangle \right|  \leq\int_{\Lambda}\left|
\left\langle f,x\right\rangle \right|  d\mu\leq\left(  \int_{\Lambda}%
Fd\mu\right)  \left|  \left|  x\right|  \right|  $, hence $\left|  \left|
\int_{\Lambda}fd\mu\right|  \right|  \leq\int_{\Lambda}Fd\mu$.

Note that if a $\ast$-dynamical system $\left(  A,\omega,\tau,K\right)  $ is
ergodic with respect to some \textit{F\o lner} net in a F\o lner semigroup
$K$, then from the mean ergodic theorem it follows that it is ergodic with
respect to every F\o lner net in $K$, since ergodicity is equivalent to the
projection $P$ in any GNS representation being $\Omega\otimes\Omega
=\Omega\left\langle \Omega,\cdot\right\rangle $, i.e. ergodicity is
independent of the F\o lner net being used, and then we can simply say that
the system is \textit{ergodic}. By the following results we have a similar
situation for weak mixing:

\bigskip

\noindent\textbf{Proposition 3.4.}\textit{ Let }$\left(  A,\omega
,\tau,K\right)  $\textit{ be a }$\ast$\textit{-dynamical system, and consider
a net }$\left(  \Lambda_{\alpha}\right)  $\textit{ in }$\Sigma$\textit{ with
}$0<\mu(\Lambda_{\alpha})<\infty$\textit{ for }$\alpha$\textit{ large enough.
Then in any GNS representation of }$\left(  A,\omega,\tau,K\right)  $\textit{
we have the following: }$\left(  A,\omega,\tau,K\right)  $\textit{ is weakly
mixing with respect to }$\left(  \Lambda_{\alpha}\right)  $\textit{ if and
only if }%
\begin{equation}
\lim_{\alpha}\frac{1}{\mu\left(  \Lambda_{\alpha}\right)  }\int_{\Lambda
_{\alpha}}\left|  \left\langle x,U_{g}y\right\rangle -\left\langle
x,\Omega\right\rangle \left\langle \Omega,y\right\rangle \right|  dg=0
\tag{3.2}%
\end{equation}
\textit{for all }$x,y\in H$\textit{.}

\bigskip

\noindent\textbf{Proof.} By setting $x=\iota\left(  a^{\ast}\right)  $ and
$y=\iota(b)$, weak mixing with respect to $\left(  \Lambda_{\alpha}\right)  $
follows immediately from (3.2).

Conversely, consider any $x,y\in H$, then there are sequences $\left(
a_{n}\right)  $ and $\left(  b_{n}\right)  $ in $A$ such that $\iota\left(
a_{n}\right)  \rightarrow x$ and $\iota\left(  b_{n}\right)  \rightarrow y$.
Hence $K\rightarrow\mathbb{C}:g\mapsto\left\langle x,U_{g}y\right\rangle
=\lim_{n\rightarrow\infty}\omega\left(  a_{n}^{\ast}\tau_{g}(b_{n})\right)  $
is $\Sigma$-measurable. Now we follow a standard argument from measure
theoretic ergodic theory (see for example \cite[Theorem 1.23]{W}). Consider
any $\varepsilon>0$, and an $n$ for which $\left|  \left|  \iota\left(
a_{n}\right)  -x\right|  \right|  <\varepsilon$ and $\left|  \left|
\iota\left(  b_{n}\right)  -y\right|  \right|  <\varepsilon$. By Definition
2.4(2) there is a $\beta$ in the directed set such that
\[
\frac{1}{\mu\left(  \Lambda_{\alpha}\right)  }\int_{\Lambda_{\alpha}}\left|
\left\langle \iota\left(  a_{n}\right)  ,U_{g}\iota\left(  b_{n}\right)
\right\rangle -\left\langle \iota\left(  a_{n}\right)  ,\Omega\right\rangle
\left\langle \Omega,\iota\left(  b_{n}\right)  \right\rangle \right|
dg<\varepsilon
\]
for every $\alpha>\beta$. By the triangle inequality
\begin{align*}
&  \left|  \left\langle x,U_{g}y\right\rangle -\left\langle x,\Omega
\right\rangle \left\langle \Omega,y\right\rangle \right| \\
&  \leq\left(  \left\|  x\right\|  +\left\|  \iota\left(  b_{n}\right)
\right\|  +\left\|  \iota\left(  a_{n}\right)  \right\|  +\left\|  y\right\|
\right)  \varepsilon+\left|  \left\langle \iota\left(  a_{n}\right)
,U_{g}\iota\left(  b_{n}\right)  \right\rangle -\left\langle \iota\left(
a_{n}\right)  ,\Omega\right\rangle \left\langle \Omega,\iota\left(
b_{n}\right)  \right\rangle \right|
\end{align*}
hence for $\alpha>\beta$ we have
\[
\frac{1}{\mu\left(  \Lambda_{\alpha}\right)  }\int_{\Lambda_{\alpha}}\left|
\left\langle x,U_{g}y\right\rangle -\left\langle x,\Omega\right\rangle
\left\langle \Omega,y\right\rangle \right|  dg\leq\left(  \left\|  x\right\|
+\left\|  \iota\left(  b_{n}\right)  \right\|  +\left\|  \iota\left(
a_{n}\right)  \right\|  +\left\|  y\right\|  +1\right)  \varepsilon
\]
while $\left\|  \iota\left(  a_{n}\right)  \right\|  \leq\left\|  x\right\|
+\varepsilon$ and $\left\|  \iota\left(  b_{n}\right)  \right\|  \leq\left\|
y\right\|  +\varepsilon$. $\square$

\bigskip

\noindent\textbf{Corollary 3.5.}\textit{ Let }$\left(  A,\omega,\tau,K\right)
$\textit{ be a }$\ast$\textit{-dynamical system with }$K$\textit{ an abelian
F\o lner semigroup, and consider any F\o lner net }$\left(  \Lambda_{\alpha
}\right)  $\textit{ in }$K$\textit{. Then in any GNS representation of
}$\left(  A,\omega,\tau,K\right)  $\textit{ we have the following: Firstly,
}$\mathbb{C}\Omega\subset H_{0}$\textit{. Secondly, }$\left(  A,\omega
,\tau,K\right)  $\textit{ is weakly mixing with respect to }$\left(
\Lambda_{\alpha}\right)  $\textit{ if and only if }$\dim H_{0}=1$\textit{. In
particular, if }$\left(  A,\omega,\tau,K\right)  $\textit{ is weakly mixing
with respect to some F\o lner net in }$K$\textit{, then it is weakly mixing
with respect to every F\o lner net in }$K$\textit{.}

\bigskip

\noindent\textbf{Proof.} Firstly note that $\mathbb{C}\Omega\subset H_{0}$,
since $U_{g}\Omega=\Omega$. Secondly, as in Proposition 3.4's proof
$K\rightarrow\mathbb{C}:g\mapsto\left\langle x,U_{g}y\right\rangle $ is
$\left(  \Sigma|_{K}\right)  $-measurable, so in particular for any $y\in H$
orthogonal to $\Omega$, Proposition 3.4 tells us that
\[
\lim_{\alpha}\frac{1}{\mu\left(  \Lambda_{\alpha}\right)  }\int_{\Lambda
_{\alpha}}\left|  \left\langle x,U_{g}y\right\rangle \right|  dg=0
\]
for all $x\in H$ when $\left(  A,\omega,\tau,K\right)  $ is weakly mixing.
Hence $y\in H_{v}$ by Proposition 3.2, so $H_{0}^{\bot}\subset\left(
\mathbb{C}\Omega\right)  ^{\bot}\subset H_{v}$, but $H_{0}^{\bot}=H_{v}$ by
Definition 3.1, hence $H_{0}=\mathbb{C}\Omega$. Conversely, if $H_{0}%
=\mathbb{C}\Omega$, then for any $x,y\in H$ write $x=x_{0}+x_{v}$ and
$y=y_{0}+y_{v}$ with $x_{0},y_{0}\in H_{0}$ and $x_{v},y_{v}\in H_{v}$. Then
it follows from Proposition 3.2 that
\begin{align*}
&  \lim_{\alpha}\frac{1}{\mu\left(  \Lambda_{\alpha}\right)  }\int
_{\Lambda_{\alpha}}\left|  \left\langle x,U_{g}y\right\rangle -\left\langle
x,\Omega\right\rangle \left\langle \Omega,y\right\rangle \right|  dg\\
&  =\lim_{\alpha}\frac{1}{\mu\left(  \Lambda_{\alpha}\right)  }\int
_{\Lambda_{\alpha}}\left|  \left\langle x,U_{g}y_{v}\right\rangle
+\left\langle x_{0},y_{0}\right\rangle -\left\langle x_{0},\left(
\Omega\otimes\Omega\right)  y_{0}\right\rangle \right|  dg\\
&  =0
\end{align*}
hence $\left(  A,\omega,\tau,K\right)  $ is weakly mixing by Proposition 3.4.
$\square$

\bigskip

Hence when $\left(  A,\omega,\tau,K\right)  $ is a $\ast$-dynamical system
which is weakly mixing with respect to some F\o lner net in an abelian F\o
lner semigroup $K$, we will simply call $\left(  A,\omega,\tau,K\right)  $
\textit{weakly mixing}.

\section{Preliminary limits}

The goal of this section is to prove Proposition 4.3, which can be viewed as a
collection of very simple nonconventional ergodic averages and is one of the
tools used in Section 5.

\bigskip

\noindent\textbf{Lemma 4.1.} \textit{Let }$\left(  A,\omega,\tau,K\right)
$\textit{\ be a C*-dynamical system which is ergodic with respect to some net
}$\left(  \Lambda_{\alpha}\right)  $\textit{\ in }$\Sigma|_{K}$\textit{. Then
}%
\[
\lim_{\alpha}\frac{1}{\mu\left(  \Lambda_{\alpha}^{-1}\right)  }\int
_{\Lambda_{\alpha}^{-1}}\omega\left(  a\tau_{g}(b)\right)  dg=\omega
(a)\omega(b)
\]
\textit{for all }$a,b\in A$\textit{. That is to say, }$\left(  A,\omega
,\tau,K^{-1}\right)  $\textit{\ is ergodic with respect to }$\left(
\Lambda_{\alpha}^{-1}\right)  $\textit{.}

\bigskip

\noindent\textbf{Proof.} Note that $\omega\left(  a\tau_{(\cdot)}(b)\right)
|_{K^{-1}}=\omega\left(  \tau_{\mathcal{I}(\cdot)}(a)b\right)  |_{K^{-1}%
}=\overline{\omega\left(  b^{\ast}\tau_{(\cdot)}(a^{\ast})\right)  |_{K}}%
\circ\mathcal{I}|_{K^{-1}}$ is $\left(  \Sigma|_{K^{-1}}\right)  $-measurable
by definition of a C*-dynamical system. Now we perform a simple calculation:%
\begin{align*}
\lim_{\alpha}\frac{1}{\mu\left(  \Lambda_{\alpha}^{-1}\right)  }\int
_{\Lambda_{\alpha}^{-1}}\omega\left(  a\tau_{g}(b)\right)  dg  &
=\lim_{\alpha}\frac{1}{\mu\left(  \Lambda_{\alpha}\right)  }\int
_{\mathcal{I}^{-1}\left(  \Lambda_{\alpha}\right)  }\omega\left(
\tau_{\mathcal{I}(g)}(a)b\right)  d\mu(g)\\
&  =\lim_{\alpha}\frac{1}{\mu\left(  \Lambda_{\alpha}\right)  }\int
_{\Lambda_{\alpha}}\omega\left(  \tau_{g}(a)b\right)  d\left(  \mu
\circ\mathcal{I}^{-1}\right)  (g)\\
&  =\overline{\lim_{\alpha}\frac{1}{\mu\left(  \Lambda_{\alpha}\right)  }%
\int_{\Lambda_{\alpha}}\omega\left(  b^{\ast}\tau_{g}(a^{\ast})\right)  dg}\\
&  =\overline{\omega(b^{\ast})\omega(a^{\ast})}\text{.\ }\square
\end{align*}

\bigskip

Note that in this lemma $G$ need not be abelian as in Definition 2.5, while
$A$ need only be a unital $\ast$-algebra. However, in the remaining results
this is not the case, hence Lemma 4.1 as it stands is good enough for our purposes.

\bigskip

\noindent\textbf{Lemma 4.2.} \textit{Let }$\left(  A,\omega,\tau,K\right)
$\textit{\ be an ergodic C*-dynamical system such that }%
\[
K\times K\rightarrow\mathbb{C}:(g,h)\mapsto\omega\left(  a\tau_{g}%
(b)\tau_{h^{-1}g}(c)\right)
\]
\textit{is }$\left(  \Sigma|_{K}\times\Sigma|_{K}\right)  $\textit{-measurable
for all }$a,b,c\in A$\textit{. Then for any F\o lner nets }$\left(
\Lambda_{1\alpha}\right)  $\textit{\ and }$\left(  \Lambda_{2\alpha}\right)
$\textit{, indexed by the same directed set, in }$K$\textit{\ and any GNS
representation of }$\left(  A,\omega\right)  $\textit{\ we have }%
\begin{equation}
\lim_{\alpha}\sup_{g_{1},g_{2}\in K}\left|  \left|  M_{\alpha,g_{1},g_{2}%
}\left(  \iota\left(  \tau_{g}(b)\tau_{h^{-1}g}(c)\right)  \right)
-\omega(b)\omega(c)\Omega\right|  \right|  =0 \tag{4.1}%
\end{equation}
\textit{for all }$b,c\in A$\textit{, where we write }%
\[
M_{\alpha,g_{1},g_{2}}\left(  f(g,h)\right)  \equiv\frac{1}{\mu\left(
\Lambda_{1\alpha}\right)  \mu\left(  \Lambda_{2\alpha}\right)  }\int
_{\Lambda_{1\alpha}g_{1}}\int_{\Lambda_{2\alpha}g_{2}}f(g,h)dgdh
\]
\textit{(in particular the symbol }$g$\textit{ on the left hand side indicates
the integration variable over }$\Lambda_{2\alpha}g_{2}$\textit{ and }%
$h$\textit{ the integration variable over }$\Lambda_{1\alpha}g_{1}$\textit{).}

\bigskip

\noindent\textbf{Proof.} $M_{\alpha,g_{1},g_{2}}\left(  \iota\left(  \tau
_{g}(b)\tau_{h^{-1}g}(c)\right)  \right)  $ exists for $\alpha$ large enough
by Fubini's Theorem, since $(g,h)\mapsto\omega\left(  a\tau_{g}(b)\tau
_{h^{-1}g}(c)\right)  =\left\langle \iota\left(  a^{\ast}\right)
,\iota\left(  \tau_{g}(b)\tau_{h^{-1}g}(c)\right)  \right\rangle $ provides
the required measurability.

Setting $c^{\prime}:=c-\omega(c)$ we have $\omega(c^{\prime})=0$, and assuming
(4.1) holds for $c^{\prime}$ instead of $c$, and keeping in mind the
invariance of $\mu$, we have%
\begin{align*}
&  \sup_{g_{1},g_{2}\in K}\left\|  M_{\alpha,g_{1},g_{2}}\left(  \iota\left(
\tau_{g}(b)\tau_{h^{-1}g}(c)\right)  \right)  -\omega(b)\omega(c)\Omega
\right\| \\
&  \leq\sup_{g_{1},g_{2}\in K}\left\|  M_{\alpha,g_{1},g_{2}}\left(
\iota\left(  \tau_{g}(b)\tau_{h^{-1}g}(c^{\prime})\right)  \right)  \right\|
+\left|  \omega(c)\right|  \sup_{g_{2}\in K}\left\|  \frac{1}{\mu\left(
\Lambda_{2\alpha}\right)  }\int_{\Lambda_{2\alpha}g_{2}}\iota\left(  \tau
_{g}(b)\right)  dg-\omega(b)\Omega\right\| \\
&  \rightarrow0
\end{align*}
in the $\alpha$ limit by the mean ergodic theorem, since $\left(
A,\omega,\tau,K\right)  $ is ergodic and $K$ is abelian, namely%
\begin{align*}
\sup_{g_{2}\in K}\left\|  \frac{1}{\mu\left(  \Lambda_{2\alpha}\right)  }%
\int_{\Lambda_{2\alpha}g_{2}}\iota\left(  \tau_{g}(b)\right)  dg-\omega
(b)\Omega\right\|   &  =\sup_{g_{2}\in K}\left\|  \frac{1}{\mu\left(
\Lambda_{2\alpha}\right)  }\int_{\Lambda_{2\alpha}}\iota\left(  \tau_{gg_{2}%
}(b)\right)  dg-\omega(b)\Omega\right\| \\
&  \leq\sup_{g_{2}\in K}\left\|  U_{g_{2}}\right\|  \left\|  \frac{1}%
{\mu\left(  \Lambda_{2\alpha}\right)  }\int_{\Lambda_{2\alpha}}U_{g}%
\iota\left(  b\right)  dg-\omega(b)\Omega\right\| \\
&  \leq\left\|  \frac{1}{\mu\left(  \Lambda_{2\alpha}\right)  }\int
_{\Lambda_{2\alpha}}U_{g}\iota\left(  b\right)  dg-\omega(b)\Omega\right\| \\
&  \rightarrow\left\|  \left(  \Omega\otimes\Omega\right)  \iota
(b)-\omega(b)\Omega\right\| \\
&  =0
\end{align*}
hence without loss we can assume $\omega(c)=0$. From our hypothesis
$h\mapsto\omega\left(  a\tau_{h^{-1}}(c)\right)  =\omega\left(  a\tau
_{e}(1)\tau_{h^{-1}e}(c)\right)  $ is $\left(  \Sigma_{K}\right)
$-measurable. Also, since $A$ is a C*-algebra and $\pi$ a $\ast$-homomorphism,
we have $\pi(a)\in B(H)$ and $\left|  \left|  \pi(a)\right|  \right|
\leq\left|  \left|  a\right|  \right|  $ for all $a\in A$. Therefore
\begin{align*}
&  \sup_{g_{1},g_{2}\in K}\left\|  M_{\alpha,g_{1},g_{2}}\left(  \iota\left(
\tau_{g}(b)\tau_{h^{-1}g}(c)\right)  \right)  \right\| \\
&  =\sup_{g_{1},g_{2}\in K}\left\|  \frac{1}{\mu\left(  \Lambda_{1\alpha
}\right)  \mu\left(  \Lambda_{2\alpha}\right)  }\int_{\Lambda_{2\alpha}g_{2}%
}\left[  \pi\left(  \tau_{g}(b)\right)  U_{g}\int_{\Lambda_{1\alpha}g_{1}%
}\iota\left(  \tau_{h^{-1}}(c)\right)  dh\right]  dg\right\| \\
&  \leq\sup_{g_{1},g_{2}\in K}\frac{1}{\mu\left(  \Lambda_{1\alpha}\right)
\mu\left(  \Lambda_{2\alpha}\right)  }\mu\left(  \Lambda_{2\alpha}%
g_{2}\right)  \left\|  b\right\|  \left\|  \int_{\Lambda_{1\alpha}g_{1}}%
\iota\left(  \tau_{h^{-1}}(c)\right)  dh\right\| \\
&  =\left|  \left|  b\right|  \right|  \sup_{g_{1}\in K}\frac{1}{\mu\left(
\Lambda_{1\alpha}\right)  }\left\|  \int_{\Lambda_{1\alpha}}\iota\left(
\tau_{\left(  hg_{1}\right)  ^{-1}}(c)\right)  dh\right\| \\
&  =\left|  \left|  b\right|  \right|  \sup_{g_{1}\in K}\frac{1}{\mu\left(
\Lambda_{1\alpha}\right)  }\left\|  U_{g_{1}^{-1}}\int_{\Lambda_{1\alpha}%
}\iota\left(  \tau_{h^{-1}}(c)\right)  dh\right\| \\
&  \leq\left|  \left|  b\right|  \right|  \frac{1}{\mu\left(  \Lambda
_{1\alpha}^{-1}\right)  }\left\|  \int_{\Lambda_{1\alpha}^{-1}}U_{h}%
\iota\left(  c\right)  dh\right\| \\
&  \rightarrow0
\end{align*}
by Lemma 4.1 and the mean ergodic theorem, since $G$ is abelian while
$K^{-1}\rightarrow\mathbb{C}:h\mapsto\left\langle \iota(a),U_{h}%
\iota(c)\right\rangle =\omega\left(  a^{\ast}\tau_{h}(c)\right)  $ is
measurable (as in Lemma 4.1's proof), ensuring the measurability required in
the last integral. $\square$

\bigskip

\noindent\textbf{Proposition 4.3.} \textit{Let }$\left(  A,\omega
,\tau,K\right)  $\textit{\ be an ergodic C*-dynamical system. Consider any F\o
lner nets }$\left(  \Lambda_{1\alpha}\right)  $\textit{\ and }$\left(
\Lambda_{2\alpha}\right)  $\textit{, over the same directed set, in }%
$K$\textit{, and use the notation }$M_{\alpha,g_{1},g_{2}}$\textit{ as in
Lemma 4.2.}

(1) \textit{If }%
\[
K\times K\rightarrow\mathbb{C}:(g,h)\mapsto\omega\left(  \tau_{h}(a)\tau
_{gh}(b)\tau_{g}(c)\right)
\]
\textit{is }$\left(  \Sigma|_{K}\times\Sigma|_{K}\right)  $\textit{-measurable
for all }$a,b,c\in A$\textit{,\ then\ we have}%
\[
\lim_{\alpha}\sup_{g_{1},g_{2}\in K}\left|  M_{\alpha,g_{1},g_{2}}\left(
\omega\left(  \tau_{g}(a)\tau_{gh}(b)\tau_{h}(c)\right)  \right)
-\omega(a)\omega(b)\omega(c)\right|  =0
\]
\textit{for all }$a,b,c\in A$\textit{.}

(2) \textit{If }$\left(  A,\omega,\tau,K\right)  $\textit{\ is asymptotically
abelian with respect to }$\left(  \Lambda_{1\alpha}\right)  $\textit{,
}$\omega$\textit{\ is tracial, and }%
\[
K\times K\rightarrow\mathbb{C}:(g,h)\mapsto\omega\left(  \tau_{h}(a_{1}%
)\tau_{g}(a_{2})\tau_{gh}(a_{3})\tau_{g}(a_{4})\tau_{h}(a_{5})\right)
\]
\textit{is }$\left(  \Sigma|_{K}\times\Sigma|_{K}\right)  $\textit{-measurable
for all }$a_{1},a_{2},a_{3},a_{4},a_{5}\in A$\textit{, then we have}%
\[
\lim_{\alpha}M_{\alpha,e,e}\left(  \omega\left(  \tau_{h}(a_{1})\tau_{g}%
(a_{2})\tau_{gh}(a_{3})\tau_{g}(a_{4})\tau_{h}(a_{5})\right)  \right)
=\omega(a_{5}a_{1})\omega(a_{3})\omega(a_{2}a_{4})
\]
\textit{for all }$a_{1},a_{2},a_{3},a_{4},a_{5}\in A$\textit{.}

(3) \textit{If }$\left(  A,\omega,\tau,K\right)  $\textit{\ is uniformly
asymptotically abelian with respect to }$\left(  \Lambda_{1\alpha}\right)
$\textit{, }$\omega$\textit{\ is tracial, and }%
\[
K\times K\rightarrow\mathbb{C}:(g,h)\mapsto\omega\left(  \tau_{h}(a_{1}%
)\tau_{g}(a_{2})\tau_{gh}(a_{3})\tau_{g}(a_{4})\tau_{h}(a_{5})\right)
\]
\textit{is }$\left(  \Sigma|_{K}\times\Sigma|_{K}\right)  $\textit{-measurable
for all }$a_{1},a_{2},a_{3},a_{4},a_{5}\in A$\textit{, then we have}%
\[
\lim_{\alpha}\sup_{g_{1},g_{2}\in K}\left|  M_{\alpha,g_{1},g_{2}}\left(
\omega\left(  \tau_{h}(a_{1})\tau_{g}(a_{2})\tau_{gh}(a_{3})\tau_{g}%
(a_{4})\tau_{h}(a_{5})\right)  \right)  -\omega(a_{5}a_{1})\omega(a_{3}%
)\omega(a_{2}a_{4})\right|  =0
\]
\textit{for all }$a_{1},a_{2},a_{3},a_{4},a_{5}\in A$\textit{.}

\bigskip

\noindent\textbf{Proof.} (1) Note that $\omega\left(  \tau_{h}(a)\tau
_{gh}(b)\tau_{g}(c)\right)  =\omega\left(  a\tau_{g}(b)\tau_{h^{-1}%
g}(c)\right)  $ therefore the latter has the measurability required in Lemma
4.2, hence
\begin{align*}
&  \sup_{g_{1},g_{2}\in K}\left|  M_{\alpha,g_{1},g_{2}}\left(  \omega\left(
\tau_{h}(a)\tau_{gh}(b)\tau_{g}(c)\right)  \right)  -\omega(a)\omega
(b)\omega(c)\right| \\
&  =\sup_{g_{1},g_{2}\in K}\left|  \left\langle \iota(a^{\ast}),M_{\alpha
,g_{1},g_{2}}\left(  \iota\left(  \tau_{g}(b)\tau_{h^{-1}g}(c)\right)
\right)  -\omega(b)\omega(c)\Omega\right\rangle \right| \\
&  \leq\left\|  \iota(a^{\ast})\right\|  \sup_{g_{1},g_{2}\in K}\left\|
M_{\alpha,g_{1},g_{2}}\left(  \iota\left(  \tau_{g}(b)\tau_{h^{-1}%
g}(c)\right)  \right)  -\omega(b)\omega(c)\Omega\right\| \\
&  \rightarrow0\text{.}%
\end{align*}

(3) Note that $(g,h)\mapsto\omega\left(  \tau_{h}(a)\tau_{g}(1)\tau
_{gh}(b)\tau_{g}(c)\tau_{h}(1)\right)  =$ $\omega\left(  a\tau_{g}%
(b)\tau_{h^{-1}g}(c)\right)  $ is $\left(  \Sigma|_{K}\times\Sigma
|_{K}\right)  $-measurable by hypothesis, hence we can apply Lemma 4.2 to
obtain
\begin{align*}
&  \sup_{g_{1},g_{2}\in K}\left|  M_{\alpha,g_{1},g_{2}}\left(  \omega\left(
\tau_{h}(a_{1})\tau_{g}(a_{2})\tau_{gh}(a_{3})\tau_{g}(a_{4})\tau_{h}%
(a_{5})\right)  \right)  -\omega(a_{5}a_{1})\omega(a_{3})\omega(a_{2}%
a_{4})\right| \\
&  =\sup_{g_{1},g_{2}\in K}\left|  \left\langle \iota\left(  \left(
a_{5}a_{1}\right)  ^{\ast}\right)  ,M_{\alpha,g_{1},g_{2}}\left(  \iota\left(
\tau_{h^{-1}g}(a_{2})\tau_{g}(a_{3})\tau_{h^{-1}g}(a_{4})\right)  \right)
-\omega(a_{3})\omega(a_{2}a_{4})\Omega\right\rangle \right| \\
&  \leq\left\|  \iota\left(  \left(  a_{5}a_{1}\right)  ^{\ast}\right)
\right\|  \sup_{g_{1},g_{2}\in K}\left\|  M_{\alpha,g_{1},g_{2}}\left(
\iota\left(  \tau_{h^{-1}g}(a_{2})\tau_{g}(a_{3})\tau_{h^{-1}g}(a_{4}%
)-\tau_{g}(a_{3})\tau_{h^{-1}g}(a_{2}a_{4})\right)  \right)  \right\| \\
&  +\left\|  \iota\left(  \left(  a_{5}a_{1}\right)  ^{\ast}\right)  \right\|
\sup_{g_{1},g_{2}\in K}\left\|  M_{\alpha,g_{1},g_{2}}\left(  \iota\left(
\tau_{g}(a_{3})\tau_{h^{-1}g}(a_{2}a_{4})\right)  \right)  -\omega
(a_{3})\omega(a_{2}a_{4})\Omega\right\| \\
&  \rightarrow0\text{,}%
\end{align*}
since
\begin{align*}
&  \left\|  \iota\left(  \tau_{h^{-1}g}(a_{2})\tau_{g}(a_{3})\tau_{h^{-1}%
g}(a_{4})-\tau_{g}(a_{3})\tau_{h^{-1}g}(a_{2}a_{4})\right)  \right\| \\
&  =\left\|  \iota\left(  \tau_{h^{-1}g}\left(  \left[  a_{2},\tau_{h}%
(a_{3})\right]  \right)  \tau_{h^{-1}g}(a_{4})\right)  \right\| \\
&  \leq\left\|  \left[  a_{2},\tau_{h}(a_{3})\right]  \right\|  \left\|
a_{4}\right\|
\end{align*}
and then by uniform asymptotic abelianness with respect to $\left(
\Lambda_{1\alpha}\right)  $%
\begin{align*}
&  \sup_{g_{1},g_{2}\in K}\left\|  M_{\alpha,g_{1},g_{2}}\left(  \iota\left(
\tau_{h^{-1}g}(a_{2})\tau_{g}(a_{3})\tau_{h^{-1}g}(a_{4})-\tau_{g}(a_{3}%
)\tau_{h^{-1}g}(a_{2}a_{4})\right)  \right)  \right\| \\
&  \leq\left\|  a_{4}\right\|  \sup_{g_{1}\in K}\frac{1}{\mu\left(
\Lambda_{1\alpha}\right)  }\int_{\Lambda_{1\alpha}g_{1}}\left\|  \left[
a_{2},\tau_{h}(a_{3})\right]  \right\|  dh\\
&  \rightarrow0\text{.}%
\end{align*}

(2) As for (3), but without the $\sup$ 's. $\square$

\section{Main results}

Now we will need a bit more structure in our F\o lner nets, namely we will
consider F\o lner nets $\left(  \Lambda_{\alpha}\right)  $ such that $\left(
\Lambda_{\alpha}^{-1}\Lambda_{\alpha}\right)  $ is also F\o lner (in such
cases $\Lambda_{\alpha}^{-1}\Lambda_{\alpha}$ need only be measurable for
$\alpha$ large enough). A simple example of this is $\left(  \Lambda
_{r}\right)  _{r\in\mathbb{R}^{+}}$, where $\Lambda_{r}:=\left\{
g\in\mathbb{R}^{q}:\left|  \left|  g\right|  \right|  <r\right\}  $, in the
group $\mathbb{R}^{q}$, since $\Lambda_{r}^{-1}\Lambda_{r}=\Lambda_{2r}$,
which also implies that $\left(  \Lambda_{r}\right)  $ satisfies the
Tempel'man condition. Furthermore, in this example both $\left(  \Lambda
_{r}\right)  $ and $\left(  \Lambda_{r}^{-1}\Lambda_{r}\right)  $ are
uniformly F\o lner in $\mathbb{R}^{q}$. The same is true if we replace
$\Lambda_{r}$ by its closure, and clearly also when $r$ is restricted to
$\mathbb{N}$ to give a sequence. Another simple example with these properties
is $\Lambda_{n}:=\left\{  1,...,n\right\}  $ in $\mathbb{Z}$.

We also need to work with second countable topological groups, i.e. groups
with a topology which has a countable topological base and in which the
product and inverse operations are continuous, because of the following ``van
der Corput'' lemma:

\bigskip

\noindent\textbf{Lemma 5.1.}\textit{ Let }$G^{\prime}$\textit{ be a second
countable topological group with a right invariant measure }$\mu^{\prime}%
$\textit{ on its }$\sigma$\textit{-algebra }$\Sigma^{\prime}$ \textit{of Borel
sets. Assume the existence of a uniform F\o lner sequence }$\left(
\Lambda_{n}^{\prime}\right)  $ \textit{in a Borel measurable semigroup
}$K^{\prime}\subset G^{\prime}$\textit{, let }$H$\textit{ be a Hilbert space,
and let }$f:G^{\prime}\rightarrow H$\textit{ be a bounded function with
}$\left\langle f(\cdot),x\right\rangle $\textit{ and }$\left\langle
f(\cdot),f(\cdot)\right\rangle :G^{\prime}\times G^{\prime}\rightarrow
\mathbb{C}$\textit{ Borel measurable for all }$x\in H$\textit{. Assume that }%
\[
\gamma_{h}:=\lim_{n\rightarrow\infty}\frac{1}{\mu^{\prime}\left(  \Lambda
_{n}^{\prime}\right)  }\int_{\Lambda_{n}^{\prime}}\left\langle
f(g),f(gh)\right\rangle dg
\]
\textit{exists for all }$h\in G$\textit{. Also assume that we have the
following limit (the given iterated integral automatically exists for }%
$m$\textit{ large enough, by the other assumptions)}%
\begin{equation}
\lim_{m\rightarrow\infty}\frac{1}{\mu(\Lambda_{m}^{\prime})^{2}}\int
_{\Lambda_{m}^{\prime}}\int_{\Lambda_{m}^{\prime}}\gamma_{h_{1}^{-1}h_{2}%
}dh_{1}dh_{2}=0\text{\textit{.}} \tag{5.1}%
\end{equation}

(1) \textit{Then it follows that }%
\[
\lim_{n\rightarrow\infty}\frac{1}{\mu^{\prime}\left(  \Lambda_{n}^{\prime
}\right)  }\int_{\Lambda_{n}^{\prime}}fd\mu^{\prime}=0\text{.}%
\]

(2) \textit{Assume furthermore that }$K^{\prime}$\textit{ is abelian, }%
$\mu^{\prime}$\textit{ is the counting measure (i.e. }$\mu^{\prime}(\Lambda
)$\textit{ is the number of elements in the set }$\Lambda\in\Sigma^{\prime}%
$\textit{, which means that }$\Lambda_{n}^{\prime}$\textit{ is a finite
non-empty set for }$n$\textit{ large enough) and that }%
\[
\lim_{n\rightarrow\infty}\sup_{g_{1}\in K^{\prime}}\left|  \frac{1}%
{\mu^{\prime}\left(  \Lambda_{n}^{\prime}\right)  }\int_{\Lambda_{n}^{\prime
}g_{1}}\left\langle f(g),f(gh)\right\rangle dg-\gamma_{h}\right|  =0
\]
\textit{for all }$h\in G$\textit{. Then}
\[
\lim_{n\rightarrow\infty}\sup_{g_{1}\in K^{\prime}}\left|  \left|
\frac{1}{\mu^{\prime}\left(  \Lambda_{n}^{\prime}\right)  }\int_{\Lambda
_{n}^{\prime}g_{1}}fd\mu^{\prime}\right|  \right|  =0\text{\textit{.}}%
\]
\textit{(The integrals in (2) are in fact finite sums, but for consistent
notation we use integral signs.)}

\bigskip

\noindent\textbf{Proof.} Part (1) was proved in \cite{BDS} following the basic
structure of a proof of a special case given in \cite{F90}. We now give an
outline of this proof, but at appropriate points we also show the minor
modifications needed to prove (2). We present it in several steps, with more
details of each step to be found in \cite{BDS}. Some steps work for nets
rather than just sequences, and in such steps we use the notation
$\Lambda_{\alpha}^{\prime}$, $\Lambda_{\beta}^{\prime}$ rather than
$\Lambda_{m}^{\prime}$, $\Lambda_{n}^{\prime}$. The fact that $G^{\prime}$ is
a topological group helps to ensure Borel measurability of the various
mappings in the proof via continuity of the group operations, while its second
countability ensures that the product $\sigma$-algebra and Borel $\sigma
$-algebra on $G^{\prime}\times G^{\prime}$ are the same, enabling us to use
Fubini's Theorem.

(a) For $\alpha$ large enough, one obtains%
\[
\lim_{\beta}\left|  \left|  \frac{1}{\mu^{\prime}(\Lambda_{\beta}^{\prime}%
)}\int_{\Lambda_{\beta}^{\prime}}fd\mu^{\prime}-\frac{1}{\mu^{\prime}%
(\Lambda_{\beta}^{\prime})}\frac{1}{\mu^{\prime}(\Lambda_{\alpha}^{\prime}%
)}\int_{\Lambda_{\beta}^{\prime}}\int_{\Lambda_{\alpha}^{\prime}%
}f(gh)dhdg\right|  \right|  =0
\]
and if furthermore $K^{\prime}$ is abelian we find that%
\[
\lim_{\beta}\sup_{g_{1}\in K}\left|  \left|  \frac{1}{\mu^{\prime}%
(\Lambda_{\beta}^{\prime})}\int_{\Lambda_{\beta}^{\prime}g_{1}}fd\mu^{\prime
}-\frac{1}{\mu^{\prime}(\Lambda_{\beta}^{\prime})}\frac{1}{\mu^{\prime
}(\Lambda_{\alpha}^{\prime})}\int_{\Lambda_{\beta}^{\prime}g_{1}}\int
_{\Lambda_{\alpha}^{\prime}}f(gh)dhdg\right|  \right|  =0\text{,}%
\]
since
\begin{align*}
&  \left|  \left|  \frac{1}{\mu^{\prime}(\Lambda_{\beta}^{\prime})}%
\int_{\Lambda_{\beta}^{\prime}g_{1}}fd\mu^{\prime}-\frac{1}{\mu^{\prime
}(\Lambda_{\beta}^{\prime})}\frac{1}{\mu^{\prime}(\Lambda_{\alpha}^{\prime}%
)}\int_{\Lambda_{\beta}^{\prime}g_{1}}\int_{\Lambda_{\alpha}^{\prime}%
}f(gh)dhdg\right|  \right| \\
&  \leq\frac{b_{1}}{\mu^{\prime}\left(  \Lambda_{\beta}^{\prime}\right)  }%
\sup_{h\in\Lambda_{\alpha}^{\prime}}\mu^{\prime}\left(  \left(  \Lambda
_{\beta}^{\prime}g_{1}\right)  \Delta\left(  \Lambda_{\beta}^{\prime}%
g_{1}h\right)  \right) \\
&  =\frac{b_{1}}{\mu^{\prime}\left(  \Lambda_{\beta}^{\prime}\right)  }%
\sup_{h\in\Lambda_{\alpha}^{\prime}}\mu^{\prime}\left(  \left(  \Lambda
_{\beta}^{\prime}\right)  \Delta\left(  \Lambda_{\beta}^{\prime}h\right)
\right)
\end{align*}
where $b_{1}$ is an upper bound for $\left|  \left|  f(K^{\prime})\right|
\right|  $ .

(b) For any $\Lambda_{1},\Lambda_{2}\in\Sigma^{\prime}|_{K^{\prime}}$ with
$\mu^{\prime}\left(  \Lambda_{1}\right)  ,\mu^{\prime}\left(  \Lambda
_{2}\right)  <\infty$ we have%
\[
\left|  \left|  \int_{\Lambda_{2}}\int_{\Lambda_{1}}f(gh)dhdg\right|  \right|
^{2}\leq\mu^{\prime}(\Lambda_{2})\int_{\Lambda_{1}}\int_{\Lambda_{1}}%
\int_{\Lambda_{2}}\left\langle f(gh_{1}),f(gh_{2})\right\rangle dgdh_{1}dh_{2}%
\]
and in particular these iterated integrals exist.

(c) We also have%
\begin{equation}
\lim_{\beta}\frac{1}{\mu^{\prime}(\Lambda_{\beta}^{\prime})}\int
_{\Lambda_{\beta}^{\prime}}\left\langle f(gh_{1}),f(gh_{2})\right\rangle
dg=\gamma_{h_{1}^{-1}h_{2}} \tag{5.2}%
\end{equation}
for all $h_{1}\in K$ and $h_{2}\in G$. However, if
\[
\lim_{\beta}\sup_{g_{1}\in K^{\prime}}\left|  \frac{1}{\mu^{\prime}\left(
\Lambda_{\beta}^{\prime}\right)  }\int_{\Lambda_{\beta}^{\prime}g_{1}%
}\left\langle f(g),f(gh)\right\rangle dg-\gamma_{h}\right|  =0
\]
for all $h\in G$ and $K^{\prime}$ is abelian, then we in fact obtain%
\begin{equation}
\lim_{\beta}\sup_{g_{1}\in K^{\prime}}\left|  \frac{1}{\mu^{\prime}%
(\Lambda_{\beta}^{\prime})}\int_{\Lambda_{\beta}^{\prime}g_{1}}\left\langle
f(gh_{1}),f(gh_{2})\right\rangle dg-\gamma_{h_{1}^{-1}h_{2}}\right|  =0
\tag{5.3}%
\end{equation}
for all $h_{1}\in K^{\prime}$ and $h_{2}\in G$ as follows:
\begin{align*}
&  \left|  \frac{1}{\mu^{\prime}(\Lambda_{\beta}^{\prime})}\int_{\Lambda
_{\beta}^{\prime}g_{1}}\left\langle f(gh_{1}),f(gh_{2})\right\rangle
dg-\gamma_{h_{1}^{-1}h_{2}}\right| \\
&  \leq\frac{\mu^{\prime}\left(  \left(  \Lambda_{\beta}^{\prime}g_{1}\right)
\Delta\left(  \Lambda_{\beta}^{\prime}g_{1}h_{1}\right)  \right)  }%
{\mu^{\prime}(\Lambda_{\beta}^{\prime})}b_{2}+\left|  \frac{1}{\mu^{\prime
}(\Lambda_{\beta}^{\prime})}\int_{\Lambda_{\beta}^{\prime}g_{1}}\left\langle
f(g),f(gh_{1}^{-1}h_{2})\right\rangle dg-\gamma_{h_{1}^{-1}h_{2}}\right|
\end{align*}
where $b_{2}$ is an upper bound for $(g,h_{1},h_{2})\mapsto\left|
\left\langle f(g),f(gh_{1}^{-1}h_{2})\right\rangle \right|  $. But
$\mu^{\prime}\left(  \left(  \Lambda_{\beta}^{\prime}g_{1}\right)
\Delta\left(  \Lambda_{\beta}^{\prime}g_{1}h_{1}\right)  \right)  =\mu
^{\prime}\left(  \left(  \Lambda_{\beta}^{\prime}\right)  \Delta\left(
\Lambda_{\beta}^{\prime}h_{1}\right)  \right)  $, since $K^{\prime}$ is
abelian. (Note that the uniformness of the F\o lner net was not required here.)

(d) Next, from (5.2) and Lebesgue's Dominated Convergence Theorem, it follows
that
\[
\lim_{n\rightarrow\infty}\frac{1}{\mu^{\prime}(\Lambda_{n}^{\prime})}%
\int_{\Lambda}\int_{\Lambda}\int_{\Lambda_{n}^{\prime}}\left\langle
f(gh_{1}),f(gh_{2})\right\rangle dgdh_{1}dh_{2}=\int_{\Lambda}\int_{\Lambda
}\gamma_{h_{1}^{-1}h_{2}}dh_{1}dh_{2}%
\]
for $\Lambda\in\Sigma^{\prime}|_{K^{\prime}}$ with $\mu^{\prime}\left(
\Lambda\right)  <\infty$.

In the case of (2), and assuming (5.3) in the case of a sequence, we have
instead that
\[
\lim_{n\rightarrow\infty}\sup_{g_{1}\in K^{\prime}}\left|  \frac{1}%
{\mu^{\prime}(\Lambda_{n}^{\prime})}\int_{\Lambda}\int_{\Lambda}\int
_{\Lambda_{n}^{\prime}g_{1}}\left\langle f(gh_{1}),f(gh_{2})\right\rangle
dgdh_{1}dh_{2}-\int_{\Lambda}\int_{\Lambda}\gamma_{h_{1}^{-1}h_{2}}%
dh_{1}dh_{2}\right|  =0
\]
which is proven as follows: Since here $\Lambda\times\Lambda$ is a finite set
and $\mu^{\prime}$ is the counting measure, integrability over $\Lambda
\times\Lambda$ is no problem and then we have%
\begin{align*}
&  \lim_{n\rightarrow\infty}\sup_{g_{1}\in K^{\prime}}\left|  \frac{1}%
{\mu^{\prime}(\Lambda_{n}^{\prime})}\int_{\Lambda}\int_{\Lambda}\int
_{\Lambda_{n}^{\prime}g_{1}}\left\langle f(gh_{1}),f(gh_{2})\right\rangle
dgdh_{1}dh_{2}-\int_{\Lambda}\int_{\Lambda}\gamma_{h_{1}^{-1}h_{2}}%
dh_{1}dh_{2}\right| \\
&  \leq\int_{\Lambda\times\Lambda}\lim_{n\rightarrow\infty}\sup_{g_{1}\in
K^{\prime}}\left|  \frac{1}{\mu^{\prime}(\Lambda_{n}^{\prime})}\int
_{\Lambda_{n}^{\prime}g_{1}}\left\langle f(gh_{1}),f(gh_{2})\right\rangle
dgdh_{1}dh_{2}-\gamma_{h_{1}^{-1}h_{2}}\right|  d\left(  h_{1},h_{2}\right)
\text{.}%
\end{align*}

(e) Lastly, combining (a), (b) and (d) with (5.1), we obtain the required
results. $\square$

\bigskip

As a corollary to this, we have a ``two parameter'' van der Corput lemma which
will be used in the proof of Corollary 5.4:

\bigskip

\noindent\textbf{Corollary 5.2.}\textit{ Let }$G$\textit{ be a second
countable topological group with a }$\sigma$\textit{-finite right and left
invariant measure }$\mu$\textit{ on its }$\sigma$\textit{-algebra }$\Sigma$
\textit{of Borel sets. Assume the existence of uniform F\o lner sequences
}$\left(  \Lambda_{1,n}\right)  $\textit{ and }$\left(  \Lambda_{2,n}\right)
$ \textit{in a Borel measurable semigroup }$K\subset G$\textit{, let }%
$H$\textit{ be a Hilbert space, and let }$f:G^{2}\rightarrow H$\textit{ be a
bounded function with }$\left\langle f(\cdot),x\right\rangle $\textit{ and
}$\left\langle f(\cdot),f(\cdot)\right\rangle :G^{2}\times G^{2}%
\rightarrow\mathbb{C}$\textit{ Borel measurable for all }$x\in H$\textit{.
Assume that }%
\[
\gamma_{\left(  g^{\prime},h^{\prime}\right)  }:=\lim_{n\rightarrow\infty
}\frac{1}{\mu\left(  \Lambda_{1,n}\right)  \mu\left(  \Lambda_{2,n}\right)
}\int_{\Lambda_{1,n}}\int_{\Lambda_{2,n}}\left\langle f(g,h),f(gg^{\prime
},hh^{\prime})\right\rangle dgdh
\]
\textit{exists for all }$g^{\prime},h^{\prime}\in G$\textit{, that }$G\times
G\rightarrow\mathbb{C}:\left(  g^{\prime},h^{\prime}\right)  \mapsto
\gamma_{\left(  g^{\prime},h^{\prime}\right)  }$\textit{ is Borel measurable,
and that }%
\begin{equation}
\lim_{m\rightarrow\infty}\frac{1}{\mu\left(  \Lambda_{1,m}\right)  \mu\left(
\Lambda_{2,m}\right)  }\int_{\Lambda_{1,m}^{-1}\Lambda_{1,m}}\int
_{\Lambda_{2,m}^{-1}\Lambda_{2,m}}\left|  \gamma_{\left(  g,h\right)
}\right|  dgdh=0 \tag{5.4}%
\end{equation}
\textit{where we also assume that }$\Lambda_{1,m}^{-1}\Lambda_{1,m}$\textit{
and }$\Lambda_{2,m}^{-1}\Lambda_{2,m}$\textit{ are Borel for }$m$\textit{
large enough.}

(1) \textit{Then it follows that }%
\[
\lim_{n\rightarrow\infty}\frac{1}{\mu\left(  \Lambda_{1,n}\right)  \mu\left(
\Lambda_{2,n}\right)  }\int_{\Lambda_{1,n}}\int_{\Lambda_{2,n}}%
f(g,h)dgdh=0\text{.}%
\]

(2) \textit{Assume furthermore that }$K$\textit{ is abelian, }$\mu$\textit{ is
the counting measure (in particular }$G$\textit{ is countable), and that }%
\[
\lim_{m\rightarrow\infty}\sup_{g_{1},g_{2}\in K}\left|  \frac{1}{\mu\left(
\Lambda_{1,n}\right)  \mu\left(  \Lambda_{2,n}\right)  }\int_{\Lambda
_{1,n}g_{1}}\int_{\Lambda_{2,n}g_{2}}\left\langle f(g,h),f(gg^{\prime
},hh^{\prime})\right\rangle dgdh-\gamma_{\left(  g^{\prime},h^{\prime}\right)
}\right|  =0
\]
\textit{for all }$g^{\prime},h^{\prime}\in G$\textit{. Then}
\[
\lim_{n\rightarrow\infty}\sup_{g_{1},g_{2}\in K}\left|  \left|  \frac{1}%
{\mu\left(  \Lambda_{1,n}\right)  \mu\left(  \Lambda_{2,n}\right)  }%
\int_{\Lambda_{1,n}g_{1}}\int_{\Lambda_{2,n}g_{2}}f(g,h)dgdh\right|  \right|
=0\text{\textit{.}}%
\]

\bigskip

\noindent\textbf{Proof.} We use a product as in Corollary 3.3's proof. It's
easily shown than $\Lambda_{n}^{\prime}:=\Lambda_{2,n}\times\Lambda_{1,n}$
gives a uniform F\o lner sequence in $K^{\prime}:=K\times K$ viewed as a F\o
lner semigroup in the ambient group $\left(  G^{\prime},\Sigma^{\prime}%
,\mu^{\prime}\right)  :=\left(  G\times G,\Sigma\times\Sigma,\mu\times
\mu\right)  $. Since $G$ is second countable, the product $\sigma$-algebra
$\Sigma^{\prime}$ is in fact the Borel $\sigma$-algebra of $G^{\prime}$. Note
that $\left(  \Lambda_{m}^{\prime}\right)  ^{-1}\Lambda_{m}^{\prime}=\left(
\Lambda_{2,m}^{-1}\Lambda_{2,m}\right)  \times\left(  \Lambda_{1,m}%
^{-1}\Lambda_{1,m}\right)  $, hence by Fubini's Theorem and \cite[Lemma
2.8]{BDS} we have%
\begin{align*}
&  \left|  \frac{1}{\mu^{\prime}(\Lambda_{m}^{\prime})^{2}}\int_{\Lambda
_{m}^{\prime}}\int_{\Lambda_{m}^{\prime}}\gamma_{\left(  h_{1},h_{2}\right)
^{-1}\left(  h_{3},h_{4}\right)  }d\left(  h_{1},h_{2}\right)  d\left(
h_{3},h_{4}\right)  \right| \\
&  \leq\frac{1}{\mu\left(  \Lambda_{1,m}\right)  \mu\left(  \Lambda
_{2,m}\right)  }\int_{\Lambda_{1,m}^{-1}\Lambda_{1,m}}\int_{\Lambda_{2,m}%
^{-1}\Lambda_{2,m}}\left|  \gamma_{\left(  g,h\right)  }\right|  dgdh
\end{align*}
for $m$ large enough, hence (5.4) implies (5.1) for the product group
$G^{\prime}$. Now simply apply Lemma 5.1 and Fubini's Theorem. $\square$

\bigskip

\noindent\textbf{Lemma 5.3.} \textit{Let }$\left(  A,\omega,\tau,G\right)
$\textit{ be an ergodic C*-dynamical system with }$\omega$\textit{ tracial,
}$G$\textit{ a group, }$p:G\times G\rightarrow G:\left(  g,h\right)  \mapsto
gh$\textit{ measurable, i.e. }$p^{-1}\left(  \Sigma\right)  \subset
\Sigma\times\Sigma$\textit{, and such that }%
\begin{equation}
G\times G\rightarrow\mathbb{C}:(g,h)\mapsto\omega\left(  \tau_{h}(a_{1}%
)\tau_{g}(a_{2})\tau_{gh}(a_{3})\tau_{g}(a_{4})\tau_{h}(a_{5})\right)
\tag{5.5}%
\end{equation}
\textit{is }$\Sigma\times\Sigma$\textit{-measurable for all }$a_{1}%
,a_{2},a_{3},a_{4},a_{5}\in A$\textit{. Assume the existence of F\o lner nets
}$\left(  \Lambda_{1\alpha}\right)  $\textit{ and }$\left(  \Lambda_{2\alpha
}\right)  $\textit{, indexed by the same directed set, in }$G$\textit{ such
that }$\left(  \Lambda_{1\alpha}^{-1}\Lambda_{1\alpha}\right)  $\textit{ and
}$\left(  \Lambda_{2\alpha}^{-1}\Lambda_{2\alpha}\right)  $\textit{ are also
F\o lner in }$G$\textit{. Consider any GNS representation of }$\left(
A,\omega\right)  $\textit{ and any }$a,b,c\in A$\textit{ with at least one of
}$\iota(a)$\textit{, }$\iota(b)$\textit{ or }$\iota(c^{\ast})$\textit{ in
}$H_{v}$\textit{. Set }$x_{g,h}:=\iota\left(  \tau_{gh}(a)\tau_{g}(b)\tau
_{h}(c)\right)  $\textit{. Write }%
\[
\gamma_{g^{\prime},h^{\prime}}\equiv\lim_{\alpha}\frac{1}{\mu\left(
\Lambda_{1\alpha}\right)  \mu\left(  \Lambda_{2\alpha}\right)  }\int
_{\Lambda_{1\alpha}}\int_{\Lambda_{2\alpha}}\left\langle x_{g,h}%
,x_{gg^{\prime},hh^{\prime}}\right\rangle dgdh
\]
\textit{for all }$g^{\prime},h^{\prime}\in G$\textit{. }

(1) \textit{If }$\left(  A,\omega,\tau,G\right)  $\textit{ is asymptotically
abelian with respect to }$\left(  \Lambda_{1\alpha}\right)  $\textit{, then we
have}%
\begin{equation}
\gamma_{g^{\prime},h^{\prime}}=\omega\left(  a^{\ast}\tau_{g^{\prime}%
h^{\prime}}(a)\right)  \omega\left(  b^{\ast}\tau_{g^{\prime}}(b)\right)
\omega\left(  \tau_{h^{\prime}}(c)c^{\ast}\right)  \tag{5.6}%
\end{equation}
\textit{giving a }$\Sigma\times\Sigma$\textit{-measurable mapping }$G\times
G\rightarrow\mathbb{C}:\left(  g,h\right)  \mapsto\gamma_{g,h}$\textit{ such
that}%
\[
\lim_{\alpha}\frac{1}{\mu\left(  \Lambda_{1\alpha}^{-1}\Lambda_{1\alpha
}\right)  \mu\left(  \Lambda_{2\alpha}^{-1}\Lambda_{2\alpha}\right)  }%
\int_{\Lambda_{1\alpha}^{-1}\Lambda_{1\alpha}}\int_{\Lambda_{2\alpha}%
^{-1}\Lambda_{2\alpha}}\left|  \gamma_{g,h}\right|  dgdh=0\text{\textit{.}}%
\]

(2) \textit{If }$\left(  A,\omega,\tau,G\right)  $\textit{ is uniformly
asymptotically abelian with respect to }$\left(  \Lambda_{1\alpha}\right)
$\textit{, then we have }%
\[
\lim_{\alpha}\sup_{g_{1},g_{2}\in G}\left|  \frac{1}{\mu\left(  \Lambda
_{1\alpha}\right)  \mu\left(  \Lambda_{2\alpha}\right)  }\int_{\Lambda
_{1\alpha}g_{1}}\int_{\Lambda_{2\alpha}g_{2}}\left\langle x_{g,h}%
,x_{gg^{\prime},hh^{\prime}}\right\rangle dgdh-\gamma_{g^{\prime},h^{\prime}%
}\right|  =0
\]
\textit{for all }$g^{\prime},h^{\prime}\in G$\textit{.}

\bigskip

\noindent\textbf{Proof.} (1) Since
\[
\left\langle x_{g,h},x_{gg^{\prime},hh^{\prime}}\right\rangle =\omega\left(
\tau_{h}(c^{\ast})\tau_{g}(b^{\ast})\tau_{gh}\left[  a^{\ast}\tau_{g^{\prime
}h^{\prime}}(a)\right]  \tau_{g}\left[  \tau_{g^{\prime}}(b)\right]  \tau
_{h}\left[  \tau_{h^{\prime}}(c)\right]  \right)  \text{,}%
\]
eq. (5.6) follows from Proposition 4.3(2). Note that $\left[  (g,h)\mapsto
\omega\left(  a^{\ast}\tau_{gh}(a)\right)  \right]  =\left[  \omega\left(
a^{\ast}\tau_{\left(  \cdot\right)  }(a)\right)  \circ\left(  (g,h)\mapsto
gh\right)  \right]  $ is $\Sigma\times\Sigma$-measurable, since $(g,h)\mapsto
gh$ is measurable and $\left(  A,\omega,\tau,G\right)  $ is a C*-dynamical
system. Similarly for $\left(  g,h\right)  \mapsto\omega\left(  b^{\ast}%
\tau_{g}(b)\right)  $ and $\left(  g,h\right)  \mapsto\omega\left(  \tau
_{h}(c)c^{\ast}\right)  =\overline{\omega\left(  c\tau_{h}\left(  c^{\ast
}\right)  \right)  }$. So $G\times G\rightarrow\mathbb{C}:\left(  g,h\right)
\mapsto\gamma_{g,h}$ is $\Sigma\times\Sigma$-measurable. Hence
\[
J_{\alpha}\equiv\frac{1}{\mu\left(  \Lambda_{1\alpha}^{-1}\Lambda_{1\alpha
}\right)  \mu\left(  \Lambda_{2\alpha}^{-1}\Lambda_{2\alpha}\right)  }%
\int_{\Lambda_{1\alpha}^{-1}\Lambda_{1\alpha}}\int_{\Lambda_{2\alpha}%
^{-1}\Lambda_{2\alpha}}\left|  \gamma_{g,h}\right|  dgdh\text{.}%
\]
exists by Fubini's Theorem for $\alpha$ large enough, and clearly $J_{\alpha
}\geq0$. We consider the three cases $\iota(a)\in H_{v}$, $\iota(b)\in H_{v}$
and $\iota\left(  c^{\ast}\right)  \in H_{v}$ separately:

For $\iota(a)\in H_{v}$ we have%
\begin{align*}
J_{\alpha}  &  \leq\frac{\left|  \left|  b\right|  \right|  ^{2}\left|
\left|  c\right|  \right|  ^{2}}{\mu\left(  \Lambda_{1\alpha}^{-1}%
\Lambda_{1\alpha}\right)  \mu\left(  \Lambda_{2\alpha}^{-1}\Lambda_{2\alpha
}\right)  }\int_{\Lambda_{1\alpha}^{-1}\Lambda_{1\alpha}}\int_{\Lambda
_{2\alpha}^{-1}\Lambda_{2\alpha}}\left|  \omega\left(  a^{\ast}\tau
_{gh}(a)\right)  \right|  dgdh\\
&  =\frac{\left|  \left|  b\right|  \right|  ^{2}\left|  \left|  c\right|
\right|  ^{2}}{\mu\left(  \Lambda_{1\alpha}^{-1}\Lambda_{1\alpha}\right)
\mu\left(  \Lambda_{2\alpha}^{-1}\Lambda_{2\alpha}\right)  }\int
_{\Lambda_{1\alpha}^{-1}\Lambda_{1\alpha}}\int_{\Lambda_{2\alpha}^{-1}%
\Lambda_{2\alpha}}\left|  \left\langle \iota(a),U_{gh}\iota(a)\right\rangle
\right|  dgdh\\
&  \rightarrow0
\end{align*}
in the $\alpha$ limit, according to Corollary 3.3, since $\left(
\Lambda_{1\alpha}^{-1}\Lambda_{1\alpha}\right)  $ and $\left(  \Lambda
_{2\alpha}^{-1}\Lambda_{2\alpha}\right)  $ are F\o lner in $G$, and note that
$(g,h)\mapsto\left\langle \iota(a_{1}),U_{gh}\iota(a_{2})\right\rangle
=\omega\left(  a_{1}^{\ast}\tau_{gh}(a_{2})\right)  $ is $\Sigma\times\Sigma
$-measurable by the same argument as above, hence $(g,h)\mapsto\left\langle
x,U_{gh}y\right\rangle $ is $\Sigma\times\Sigma$-measurable for all $x,y\in H$
by considering sequences in $\iota(A)$ converging to $x$ and $y$ (as in
Proposition 3.4's proof).

In much the same way for $\iota(b)\in H_{v}$ we have%
\begin{align*}
J_{\alpha}  &  \leq\frac{\left|  \left|  a\right|  \right|  ^{2}\left|
\left|  c\right|  \right|  ^{2}}{\mu\left(  \Lambda_{1\alpha}^{-1}%
\Lambda_{1\alpha}\right)  \mu\left(  \Lambda_{2\alpha}^{-1}\Lambda_{2\alpha
}\right)  }\int_{\Lambda_{1\alpha}^{-1}\Lambda_{1\alpha}}\int_{\Lambda
_{2\alpha}^{-1}\Lambda_{2\alpha}}\left|  \omega\left(  b^{\ast}\tau
_{g}(b)\right)  \right|  dgdh\\
&  =\frac{\left|  \left|  a\right|  \right|  ^{2}\left|  \left|  c\right|
\right|  ^{2}}{\mu\left(  \Lambda_{2\alpha}^{-1}\Lambda_{2\alpha}\right)
}\int_{\Lambda_{2\alpha}^{-1}\Lambda_{2\alpha}}\left|  \left\langle
\iota(b),U_{g}\iota(b)\right\rangle \right|  dg\\
&  \rightarrow0
\end{align*}
in the $\alpha$ limit according to Proposition 3.2. The case $\iota\left(
c^{\ast}\right)  \in H_{v}$ is similar to $\iota(b)\in H_{v}$, but using
$\omega\left(  \tau_{h}(c)c^{\ast}\right)  =\overline{\omega\left(  c\tau
_{h}(c^{\ast})\right)  }$.

(2) This follows directly from Proposition 4.3(3), the formula for
$\left\langle x_{g,h},x_{gg^{\prime},hh^{\prime}}\right\rangle $ given in
(1)'s proof, and (5.6), since $G$ is a group and hence uniform asymptotic
abelianness implies asymptotic abelianness with respect to $\left(
\Lambda_{1\alpha}\right)  $.\textbf{ }$\square$

\bigskip

\noindent\textbf{Corollary 5.4.}\textit{ Assume the situation in Lemma 5.3,
but assume instead of (5.5) that }%
\[
G^{2}\times G^{2}\rightarrow\mathbb{C}:\left(  g,h,j,k\right)  \mapsto
\omega\left(  \tau_{h}(a_{1})\tau_{g}(a_{2})\tau_{gh}(a_{3})\tau_{jk}%
(a_{4})\tau_{j}(a_{5})\tau_{k}(a_{6})\right)
\]
\textit{is }$\Sigma^{2}\times\Sigma^{2}$\textit{-measurable for all }%
$a_{1},a_{2},a_{3},a_{4},a_{5},a_{6}\in A$\textit{, and furthermore that }%
$G$\textit{ is a second countable topological group with }$\Sigma$\textit{ its
Borel }$\sigma$\textit{-algebra (so }$p^{-1}\left(  \Sigma\right)
\subset\Sigma\times\Sigma$\textit{ is automatically satisfied), and that }%
$\mu$\textit{ is }$\sigma$\textit{-finite. Assume the existence of uniform F\o
lner sequences }$\left(  \Lambda_{1,n}\right)  $\textit{ and }$\left(
\Lambda_{2,n}\right)  $\textit{ in }$G$\textit{ satisfying the Tempel'man
condition and such that }$\left(  \Lambda_{1,n}^{-1}\Lambda_{1,n}\right)
$\textit{ and }$\left(  \Lambda_{2,n}^{-1}\Lambda_{2,n}\right)  $\textit{ are
also F\o lner in }$G$\textit{. }

(1) \textit{If }$\left(  A,\omega,\tau,G\right)  $\textit{ is asymptotically
abelian with respect to }$\left(  \Lambda_{1n}\right)  $\textit{ then }%
\[
\lim_{n\rightarrow\infty}\frac{1}{\mu\left(  \Lambda_{1,n}\right)  \mu\left(
\Lambda_{2,n}\right)  }\int_{\Lambda_{1,n}}\int_{\Lambda_{2,n}}x_{g,h}%
dgdh=0\text{\textit{.}}%
\]

(2) \textit{If }$\mu$\textit{ is the counting measure and }$\left(
A,\omega,\tau,G\right)  $\textit{ is uniformly asymptotically abelian with
respect to }$\left(  \Lambda_{1n}\right)  $\textit{, then}%
\[
\lim_{n\rightarrow\infty}\sup_{g_{1},g_{2}\in G}\left\|  \frac{1}{\mu\left(
\Lambda_{1,n}\right)  \mu\left(  \Lambda_{2,n}\right)  }\int_{\Lambda
_{1,n}g_{1}}\int_{\Lambda_{2,n}g_{2}}x_{g,h}dgdh\right\|  =0\text{\textit{.}}%
\]

\bigskip

\noindent\textbf{Proof.} The mapping (5.5) can be expressed as
\[
\left(  g,h\right)  \mapsto\left(  g,h,g,h\right)  \mapsto\omega\left(
\tau_{h}(a_{1})\tau_{g}(a_{2})\tau_{gh}(a_{3})\tau_{gh}(1)\tau_{g}(a_{4}%
)\tau_{h}(a_{5})\right)
\]
which by the hypothesis is $\Sigma\times\Sigma$-measurable as required in
Lemma 5.3, since in any topological space $y\mapsto(y,y)$ is continuous in the
product topology. Now simply apply Corollary 5.2 and the Tempel'man condition
to $f(g,h)=x_{g,h}$, noting that $\left\|  f(g,h)\right\|  \leq\left\|
\tau_{gh}(a)\tau_{g}(b)\tau_{h}(c)\right\|  \leq\left\|  a\right\|  \left\|
b\right\|  \left\|  c\right\|  $, so $f$ is bounded, and
\[
(g,h)\mapsto\left\langle f(g,h),\iota(d)\right\rangle =\omega\left(  \tau
_{h}(c^{\ast})\tau_{g}(b^{\ast})\tau_{gh}(a^{\ast})\tau_{je}(1)\tau_{j}%
(1)\tau_{e}(d)\right)
\]
is $\Sigma^{2}$-measurable, hence so is $\left(  g,h\right)  \mapsto
\left\langle f(g,h),x\right\rangle $ for all $x\in H$, since $\iota(A)$ is
dense in $H$, while
\[
(g,h,j,k)\mapsto\left\langle f(g,h),f(j,k)\right\rangle =\omega\left(
\tau_{h}(c^{\ast})\tau_{g}(b^{\ast})\tau_{gh}(a^{\ast})\tau_{jk}(a)\tau
_{j}(b)\tau_{k}(c)\right)
\]
is $\Sigma^{2}\times\Sigma^{2}$-measurable. $\square$

\bigskip

\noindent\textbf{Proposition 5.5.}\textit{ Assume the situation in Corollary
5.4, and also that }$\left(  A,\omega,\tau,G\right)  $\textit{ is weakly
mixing, but now set }%
\[
x_{g,h}(a,b,c):=\iota\left(  \tau_{gh}(a)\tau_{g}(b)\tau_{h}(c)\right)
\]
\textit{and }%
\[
L(a,b,c):=\omega(a)\omega(b)\omega(c)
\]
\textit{ for all }$a,b,c\in A$\textit{. }

(1) \textit{If }$\left(  A,\omega,\tau,G\right)  $\textit{ is asymptotically
abelian with respect to }$\left(  \Lambda_{1n}\right)  $\textit{, then }%
\[
\lim_{n\rightarrow\infty}\frac{1}{\mu\left(  \Lambda_{1,n}\right)  \mu\left(
\Lambda_{2,n}\right)  }\int_{\Lambda_{1,n}}\int_{\Lambda_{2,n}}x_{g,h}%
(a,b,c)dgdh=L(a,b,c)\Omega
\]
\textit{for all }$a,b,c\in A$.

(2) \textit{If }$\mu$\textit{ is the counting measure, and }$\left(
A,\omega,\tau,G\right)  $\textit{ is uniformly asymptotically abelian with
respect to }$\left(  \Lambda_{1n}\right)  $\textit{, then}%
\[
\lim_{n\rightarrow\infty}\sup_{g_{1},g_{2}\in G}\left\|  \frac{1}{\mu\left(
\Lambda_{1,n}\right)  \mu\left(  \Lambda_{2,n}\right)  }\sum_{h\in
\Lambda_{1,n}g_{1}}\sum_{g\in\Lambda_{2,n}g_{2}}x_{g,h}(a,b,c)-L(a,b,c)\Omega
\right\|  =0
\]
\textit{for all }$a,b,c\in A$.

\bigskip

\noindent\textbf{Proof.} Write $a_{0}:=\omega(a)1$ and $a_{v}:=a-a_{0}$, then
by Corollary 3.5 $\iota(a_{0})=\omega(a)\Omega\in H_{0}$ and $\iota(a_{v})\in
H_{v}$, since $\left\langle \Omega,\iota\left(  a_{v}\right)  \right\rangle
=\omega\left(  1^{\ast}a_{v}\right)  =0$. Similarly for $b$ and $c$, so in
particular $x_{g,h}(a_{0},b_{0},c_{0})=L(a,b,c)\Omega$. Furthermore,
$\iota\left(  c_{0}^{\ast}\right)  =\overline{\omega(c)}\Omega\in H_{0}$ and
$\left\langle \iota\left(  c_{v}^{\ast}\right)  ,\Omega\right\rangle
=\omega\left(  c_{v}1\right)  =0$ so $\iota\left(  c_{v}^{\ast}\right)  \in
H_{v}$.

(1) By Corollary 5.4(1) we then have
\begin{align*}
&  \lim_{n\rightarrow\infty}\frac{1}{\mu\left(  \Lambda_{1,n}\right)
\mu\left(  \Lambda_{2,n}\right)  }\int_{\Lambda_{1,n}}\int_{\Lambda_{2,n}%
}x_{g,h}(a,b,c)dgdh\\
&  =\lim_{n\rightarrow\infty}\frac{1}{\mu\left(  \Lambda_{1,n}\right)
\mu\left(  \Lambda_{2,n}\right)  }\int_{\Lambda_{1,n}}\int_{\Lambda_{2,n}%
}x_{g,h}(a_{0},b_{0},c_{0})dgdh\\
&  =L(a,b,c)\lim_{n\rightarrow\infty}\frac{1}{\mu\left(  \Lambda_{1,n}\right)
\mu\left(  \Lambda_{2,n}\right)  }\int_{\Lambda_{1,n}}\int_{\Lambda_{2,n}%
}\Omega dgdh\\
&  =L(a,b,c)\Omega\text{. }%
\end{align*}

(2) Similarly by Corollary 5.4(2), but switching to summation notation and
also using the triangle inequality, we have%
\begin{align*}
&  \lim_{n\rightarrow\infty}\sup_{g_{1},g_{2}\in G}\left\|  \frac{1}%
{\mu\left(  \Lambda_{1,n}\right)  \mu\left(  \Lambda_{2,n}\right)  }\sum
_{h\in\Lambda_{1,n}g_{1}}\sum_{g\in\Lambda_{2,n}g_{2}}x_{g,h}%
(a,b,c)-L(a,b,c)\Omega\right\| \\
&  \leq\lim_{n\rightarrow\infty}\sup_{g_{1},g_{2}\in G}\left\|  \frac{1}%
{\mu\left(  \Lambda_{1,n}\right)  \mu\left(  \Lambda_{2,n}\right)  }\sum
_{h\in\Lambda_{1,n}g_{1}}\sum_{g\in\Lambda_{2,n}g_{2}}x_{g,h}(a_{0}%
,b_{0},c_{0})-L(a,b,c)\Omega\right\| \\
&  =0\text{. }\square
\end{align*}

\bigskip

This proposition in itself is interesting. Translating all the requirements to
a measure theoretic system, (1) reduces to $L^{2}$ convergence of the
following nonconventional ergodic average:%
\begin{align*}
&  \lim_{n\rightarrow\infty}\frac{1}{\mu\left(  \Lambda_{1,n}\right)
\mu\left(  \Lambda_{2,n}\right)  }\int_{\Lambda_{1,n}}\int_{\Lambda_{2,n}%
}\left(  f_{1}\circ T_{gh}\right)  \left(  f_{2}\circ T_{g}\right)  \left(
f_{3}\circ T_{h}\right)  dgdh\\
&  =\left(  \int f_{1}d\nu\right)  \left(  \int f_{2}d\nu\right)  \left(  \int
f_{3}d\nu\right)
\end{align*}
for $f_{1},f_{2},f_{3}\in L^{\infty}\left(  \nu\right)  $ where $\nu$ is a
probability measure on some measurable space and with $T_{g}$ an invertible
measure preserving transformation of this probability space, keeping in mind
that a GNS\ representation is now simply given by the set inclusion
$\iota:L^{\infty}(\nu)\rightarrow L^{2}(\nu)$. Similarly for (2).

Now we arrive at the main result of this paper:

\bigskip

\noindent\textbf{Theorem 5.6.}\textit{ Let }$\left(  A,\omega,\tau,G\right)
$\textit{ be a weakly mixing C*-dynamical system with }$\omega$\textit{
tracial, }$G$\textit{ a second countable topological group, }$\Sigma$\textit{
its Borel }$\sigma$\textit{-algebra, }$\mu$\textit{ a }$\sigma$\textit{-finite
measure, and such that }%
\[
G^{2}\times G^{2}\rightarrow\mathbb{C}:\left(  g,h,j,k\right)  \mapsto
\omega\left(  \tau_{h}(a_{1})\tau_{g}(a_{2})\tau_{gh}(a_{3})\tau_{jk}%
(a_{4})\tau_{j}(a_{5})\tau_{k}(a_{6})\right)
\]
\textit{is }$\Sigma^{2}\times\Sigma^{2}$\textit{-measurable. Assume the
existence of uniform F\o lner sequences }$\left(  \Lambda_{1,n}\right)
$\textit{ and }$\left(  \Lambda_{2,n}\right)  $\textit{ in }$G$\textit{
satisfying the Tempel'man condition and such that }$\left(  \Lambda_{1,n}%
^{-1}\Lambda_{1,n}\right)  $\textit{ and }$\left(  \Lambda_{2,n}^{-1}%
\Lambda_{2,n}\right)  $\textit{ are also F\o lner in }$G$\textit{.}

(1) \textit{If }$\left(  A,\omega,\tau,G\right)  $\textit{ is asymptotically
abelian with respect to }$\left(  \Lambda_{1n}\right)  $\textit{, then }%
\[
\lim_{n\rightarrow\infty}\frac{1}{\mu\left(  \Lambda_{1,n}\right)  \mu\left(
\Lambda_{2,n}\right)  }\int_{\Lambda_{1,n}}\int_{\Lambda_{2,n}}\omega\left(
\tau_{gh}(a)\tau_{g}(b)\tau_{h}(c)d\right)  dgdh=\omega(a)\omega
(b)\omega(c)\omega(d)
\]
\textit{for all }$a,b,c,d\in A$\textit{.}

(2) \textit{If }$\mu$\textit{ is the counting measure (in particular }%
$G$\textit{ is countable), and }$\left(  A,\omega,\tau,G\right)  $\textit{ is
uniformly asymptotically abelian with respect to }$\left(  \Lambda
_{1,n}\right)  $\textit{, then}%
\[
\lim_{n\rightarrow\infty}\sup_{g_{1},g_{2}\in G}\left|  \frac{1}{\mu\left(
\Lambda_{1,n}\right)  \mu\left(  \Lambda_{2,n}\right)  }\sum_{h\in
\Lambda_{1,n}g_{1}}\sum_{g\in\Lambda_{2,n}g_{2}}\omega\left(  \tau_{gh}%
(a)\tau_{g}(b)\tau_{h}(c)d\right)  -\omega(a)\omega(b)\omega(c)\omega
(d)\right|  =0
\]
\textit{for all }$a,b,c,d\in A$\textit{.}

\bigskip

\noindent\textbf{Proof.} This follows easily from Proposition 5.5, namely for
(1) we have\textit{ }%
\begin{align*}
&  \lim_{n\rightarrow\infty}\frac{1}{\mu\left(  \Lambda_{1,n}\right)
\mu\left(  \Lambda_{2,n}\right)  }\int_{\Lambda_{1,n}}\int_{\Lambda_{2,n}%
}\omega\left(  \tau_{gh}(a)\tau_{g}(b)\tau_{h}(c)d\right)  dgdh\\
&  =\left\langle \iota(d^{\ast}),\lim_{n\rightarrow\infty}\frac{1}{\mu\left(
\Lambda_{1,n}\right)  \mu\left(  \Lambda_{2,n}\right)  }\int_{\Lambda_{1,n}%
}\int_{\Lambda_{2,n}}x_{g,h}(a,b,c)dgdh\right\rangle \\
&  =\left\langle \iota(d^{\ast}),L(a,b,c)\Omega\right\rangle \\
&  =\omega(d)\omega(a)\omega(b)\omega(c)
\end{align*}
and similarly for (2) using
\begin{align*}
&  \left|  \frac{1}{\mu\left(  \Lambda_{1,n}\right)  \mu\left(  \Lambda
_{2,n}\right)  }\sum_{h\in\Lambda_{1,n}g_{1}}\sum_{g\in\Lambda_{2,n}g_{2}%
}\omega\left(  \tau_{gh}(a)\tau_{g}(b)\tau_{h}(c)d\right)  -\omega
(a)\omega(b)\omega(c)\omega(d)\right| \\
&  =\left|  \left\langle \iota(d^{\ast}),\frac{1}{\mu\left(  \Lambda
_{1,n}\right)  \mu\left(  \Lambda_{2,n}\right)  }\sum_{h\in\Lambda_{1,n}g_{1}%
}\sum_{g\in\Lambda_{2,n}g_{2}}x_{g,h}(a,b,c)-L(a,b,c)\Omega\right\rangle
\right| \\
&  \leq\left\|  \iota(d^{\ast})\right\|  \left\|  \frac{1}{\mu\left(
\Lambda_{1,n}\right)  \mu\left(  \Lambda_{2,n}\right)  }\sum_{h\in
\Lambda_{1,n}g_{1}}\sum_{g\in\Lambda_{2,n}g_{2}}x_{g,h}(a,b,c)-L(a,b,c)\Omega
\right\|
\end{align*}
and Proposition 5.5's notation. $\square$

\bigskip

From this we can derive the following recurrence result:

\bigskip

\noindent\textbf{Corollary 5.7.}\textit{ Consider the situation in Theorem 5.6
and let }$\varepsilon>0$\textit{ be given. Consider any }$a,b,c,d\in A$\textit{.}

(1) \textit{If }$\left(  A,\omega,\tau,G\right)  $\textit{ is asymptotically
abelian with respect to }$\left(  \Lambda_{1n}\right)  $\textit{, then there
is an }$n_{0}\in\mathbb{N}$\textit{ such that for every }$n>n_{0}$\textit{ }%
\[
\left|  \omega\left(  \tau_{gh}(a)\tau_{g}(b)\tau_{h}(c)d\right)  \right|
>\left|  \omega(a)\omega(b)\omega(c)\omega(d)\right|  -\varepsilon
\]
\textit{for some }$g\in\Lambda_{2,n}$\textit{ and some }$h\in\Lambda_{1,n}$\textit{.}

(2) \textit{If }$\mu$\textit{ is the counting measure, and }$\left(
A,\omega,\tau,G\right)  $\textit{ is uniformly asymptotically abelian with
respect to }$\left(  \Lambda_{1\alpha}\right)  $\textit{, then there is an
}$n\in\mathbb{N}$\textit{ such that for every }$g_{1},g_{2}\in G$\textit{ we
have }%
\begin{equation}
\left|  \omega\left(  \tau_{gh}(a)\tau_{g}(b)\tau_{h}(c)d\right)  \right|
>\left|  \omega(a)\omega(b)\omega(c)\omega(d)\right|  -\varepsilon\tag{5.7}%
\end{equation}
\textit{for some }$g\in\Lambda_{2,n}g_{2}$\textit{ and some }$h\in
\Lambda_{1,n}g_{1}$\textit{, i.e. the set of }$\left(  g,h\right)  $\textit{
's for which (5.7) holds is relatively dense in }$G\times G$\textit{.}

\bigskip

\noindent\textbf{Proof.} We only prove (2), since (1)'s proof is similar. By
Theorem 5.6(2) there exists an $n$ such that
\[
\frac{1}{\mu\left(  \Lambda_{1,n}\right)  \mu\left(  \Lambda_{2,n}\right)
}\sum_{h\in\Lambda_{1,n}g_{1}}\sum_{g\in\Lambda_{2,n}g_{2}}\left|
\omega\left(  \tau_{gh}(a)\tau_{g}(b)\tau_{h}(c)d\right)  \right|  >\left|
\omega(a)\omega(b)\omega(c)\omega(d)\right|  -\varepsilon
\]
for all $g_{1},g_{2}\in G$, from which the result follows. Also keep in mind
that relative denseness of a set $E$ in $G\times G$ is often defined in the
following equivalent way: $FE=G\times G$ for some finite set $F$ in $G\times
G$, in this case $F=\left(  \Lambda_{2,n}\times\Lambda_{1,n}\right)  ^{-1}$.
$\square$

\bigskip

In the case of a countable group, Corollary 5.7(2) therefore says that we
``regularly'' have recurrence. In the more general situation, Corollary 5.7(1)
isn't quite as strong, however keep in mind that the intervals $\Lambda
_{1}=\left[  0,1\right]  $, $\Lambda_{2}\in\left[  1,3\right]  $, $\Lambda
_{3}\in\left[  3,6\right]  $, ... give a uniform F\o lner sequence with the
required properties in $G=\mathbb{R}$, so in this case Corollary 5.7(1) says
that from a certain interval onward, we do get recurrence in each interval,
but with the intervals steadily growing in size. Similarly in $G=\mathbb{R}%
^{q}$, where for example we can use a sequence of balls as at the beginning of
this section, but shifted so that they don't overlap.

\bigskip

\noindent\textbf{Acknowledgment}

I thank Conrad Beyers, Anton Str\"{o}h and Johan Swart for useful conversations.

\end{document}